\documentclass{article}

\newcommand{\bA}{{\mathbf{A}}}

\newcommand{\bB}{{\mathbf{B}}}

\newcommand{\bC}{{\mathbf{C}}}

\newcommand{\cD}{{\mathcal{D}}}

\newcommand{\cE}{{\mathcal{E}}}

\newcommand{\bF}{{\mathbf{F}}}
\newcommand{\cF}{{\mathcal{F}}}

\newcommand{\bG}{{\mathbf{G}}}
\newcommand{\cG}{{\mathcal{G}}}

\newcommand{\bH}{{\mathbf{H}}}

\newcommand{\bI}{{\mathbf{I}}}

\newcommand{\cN}{{\mathcal{N}}}

\newcommand{\cO}{{\mathcal{O}}}

\newcommand{\bQ}{{\mathbf{Q}}}

\newcommand{\bR}{{\mathbf{R}}}
\newcommand{\cR}{{\mathcal{R}}}

\newcommand{\bu}{{\mathbf{u}}}

\newcommand{\cU}{{\mathcal{U}}}

\newcommand{\bv}{{\mathbf{v}}}

\newcommand{\bW}{{\mathbf{W}}}

\newcommand{\bX}{{\mathbf{X}}}

\newcommand{\bY}{{\mathbf{Y}}}

\newcommand{\bpi}{{\boldsymbol{\pi}}}

\newcommand{\bxi}{{\boldsymbol{\xi}}}

\newcommand{\bPi}{{\boldsymbol{\Pi}}}

\newcommand{\expect}{\mathbb{E}}
\newcommand{\reals}{\mathbb{R}}

\newcommand{\ip}[2]{\left\langle{#1},{#2}\right\rangle}

\newcommand{\T}{{\mathrm{T}}}

\newcommand{\btheta}{\mathbf{\Theta}}
\newcommand{\tbR}{\tilde{\bR}}
\newcommand{\tbC}{\tilde{\bC}}
\newcommand{\tbW}{\tilde{\bW}}
\newcommand{\bbF}{\bar{\bF}}

\newcommand{\prt}[1]{\left(#1\right)}
\newcommand{\brk}[1]{\left[#1\right]}
\newcommand{\crk}[1]{\left\{#1\right\}}

\newcommand{\norm}[1]{\left\|#1\right\|}

\usepackage{bbm}
\newcommand{\mone}{\mathbbm{1}}
\usepackage{arxiv}
\usepackage[applemac]{inputenc} 		
\usepackage[T1]{fontenc}    		
\usepackage[colorlinks]{hyperref}      
\usepackage{url}            			
\usepackage{amsthm}
\usepackage{booktabs}       		
\usepackage{amsfonts}       		
\usepackage{amsmath}
\usepackage{nicefrac}       		
\usepackage{microtype}      		
\usepackage{lipsum}				
\usepackage[square,numbers]{natbib}
\usepackage{mathtools}
\usepackage{algpseudocode}
\usepackage{graphicx}
\usepackage{indentfirst,latexsym,bm}
\usepackage{amsmath}
\usepackage{amssymb}
\usepackage{xcolor}
\usepackage{comment}
\usepackage{enumitem}
\usepackage{dsfont}
\usepackage{amssymb}

\usepackage{algorithm}
\usepackage{algorithmicx}
\usepackage{float}
\usepackage{subfig}
\usepackage{caption}
\usepackage{enumitem}
\usepackage{todonotes}
\usepackage{makecell}
\usepackage[flushleft]{threeparttable} 

\usepackage{tikz}
\usetikzlibrary{arrows.meta,positioning}
\usepackage{pgfplots}
\pgfplotsset{compat=newest}
\usepgfplotslibrary{fillbetween}
\usetikzlibrary{shapes,decorations}
\usetikzlibrary{fit}

\colorlet{color1}{blue}
\colorlet{color2}{red!50!black}

\definecolor{ivory}{RGB}{218,215,203}

\definecolor{cuhkp}{RGB}{98,56,105} 	
\definecolor{cuhkpl}{RGB}{152,24,147} 	
\definecolor{cuhkb}{RGB}{219,160,1} 	
\definecolor{cuhkbd}{RGB}{178,129,0} 	
\definecolor{cuhkr}{RGB}{88,35,155}  	

\usepackage{hyperref}[6.83]
\hypersetup{
  colorlinks=true,
  frenchlinks=false,
  pdfborder={0 0 0},
  naturalnames=false,
  hypertexnames=false,
  breaklinks,
  allcolors = cuhkbd,
  urlcolor = color2,
}

\RequirePackage[capitalize,nameinlink]{cleveref}[0.19]

\crefname{section}{section}{sections}
\crefname{subsection}{subsection}{subsections}

\Crefname{figure}{Figure}{Figures}

\crefformat{equation}{\textup{#2(#1)#3}}
\crefrangeformat{equation}{\textup{#3(#1)#4--#5(#2)#6}}
\crefmultiformat{equation}{\textup{#2(#1)#3}}{ and \textup{#2(#1)#3}}
{, \textup{#2(#1)#3}}{, and \textup{#2(#1)#3}}
\crefrangemultiformat{equation}{\textup{#3(#1)#4--#5(#2)#6}}%
{ and \textup{#3(#1)#4--#5(#2)#6}}{, \textup{#3(#1)#4--#5(#2)#6}}{, and \textup{#3(#1)#4--#5(#2)#6}}

\Crefformat{equation}{#2Equation~\textup{(#1)}#3}
\Crefrangeformat{equation}{Equations~\textup{#3(#1)#4--#5(#2)#6}}
\Crefmultiformat{equation}{Equations~\textup{#2(#1)#3}}{ and \textup{#2(#1)#3}}
{, \textup{#2(#1)#3}}{, and \textup{#2(#1)#3}}
\Crefrangemultiformat{equation}{Equations~\textup{#3(#1)#4--#5(#2)#6}}%
{ and \textup{#3(#1)#4--#5(#2)#6}}{, \textup{#3(#1)#4--#5(#2)#6}}{, and \textup{#3(#1)#4--#5(#2)#6}}

\crefdefaultlabelformat{#2\textup{#1}#3}

\theoremstyle{plain}
\newtheorem{theorem}{Theorem}[section]
\newtheorem{thm}{Theorem}[section]
\newtheorem{lemma}[thm]{Lemma}
\newtheorem{corollary}[thm]{Corollary}
\newtheorem{definition}[thm]{Definition}

\newtheorem{remark}{Remark}[section]
\newtheorem{assumption}[thm]{Assumption}
\theoremstyle{plain}

\title{Stochastic Push-Pull for Decentralized Nonconvex Optimization}

\author{%
  Runze You\\
  School of Data Science\\
  The Chinese University of Hong Kong, Shenzhen (CUHK-Shenzhen)\\
  \texttt{runzeyou@link.cuhk.edu.cn} \\
  \And
  Shi Pu \\
  School of Data Science \\
  The Chinese University of Hong Kong, Shenzhen (CUHK-Shenzhen)\\
  \texttt{pushi@cuhk.edu.cn} \\
}
\date{}

\begin{document}
\maketitle

\begin{abstract}
    To understand the convergence behavior of the Push-Pull method for decentralized optimization with stochastic gradients (Stochastic Push-Pull), this paper presents a comprehensive analysis. Specifically, we first clarify the algorithm's underlying assumptions, particularly those regarding the network structure and weight matrices. Then, to establish the convergence rate under smooth nonconvex objectives, we introduce a general analytical framework that not only encompasses a broad class of decentralized optimization algorithms, but also recovers or enhances several state-of-the-art results for distributed stochastic gradient tracking methods. A key highlight is the derivation of a sufficient condition under which the Stochastic Push-Pull algorithm achieves linear speedup, matching the scalability of centralized stochastic gradient methods---a result not previously reported. Extensive numerical experiments validate our theoretical findings, demonstrating the algorithm's effectiveness and robustness across various decentralized optimization scenarios.
\end{abstract}

\section{Introduction}

Consider a set of agents/nodes $\cN = \{1,2,\cdots,n\}$ connected over a network. Each node has its own local cost function $f_i:\reals^p \rightarrow \reals$. The global objective is to locate $x\in \reals^p$ that minimizes the average of all cost functions, i.e.,
\begin{equation}
    \label{eq:obj}
    \min_{x\in \reals^p} f(x) := \frac{1}{n}\sum_{i=1}^n f_i(x),
\end{equation}
where $f_i(x) := \expect_{\xi_i\sim \cD_i} F_i(x;\xi_i)$. The function $f_i(x)$ is kept at node $i$, and $\cD_i$ denotes the distribution of sample $\xi_i$ locally stored at node $i$. Data heterogeneity exists if the local distributions $\crk{\cD_i}_{i=1}^{n}$ are not identical. The above formulation also covers empirical risk minimization (ERM), with each $\cD_i$ representing a local training dataset. 

Solving Problem \eqref{eq:obj} in a decentralized manner over a multi-agent network has gained great interests in recent years, partly because decentralization helps relieve the high latency at the central server in a centralized computing architecture \cite{nedic2018network,assran2019stochastic,nadiradze2021asynchronous}. In a decentralized setup, the agents form a specific network topology (e.g., ring, grid, exponential graph) and exchange information locally with only their directed neighbors \cite{nedic2009distributed,chen2012diffusion,nedic2018network}. Such communication patterns greatly reduce the communication overhead. They also accommodate the scenarios in which certain agents can only reach a subset of the network (due to varying power ranges or physical restrictions) \cite{yang2019survey}.

The critical factor distinguishing decentralized optimization methods from centralized ones is the communication graph (or the network topology). The structure of the graph significantly affects the convergence rate: denser topologies, such as exponential graphs, typically offer faster convergence compared to sparser structures, like ring graphs \cite{pu2021sharp,nedic2018network}. Notably, under stochastic gradient settings, several decentralized algorithms have been shown to eventually match the convergence rate of centralized methods (linear-speedup), but only after enduring a transient phase \cite{pu2020asymptotic,pu2021sharp,lian2017can,tang2018d,ying2021exponential}. The number of transient iterations is heavily influenced by the underlying network topology, making them a crucial measure of the algorithmic efficiency.

A common analytical tool for assessing the impact of network topology is the spectral gap of the associated weight matrix, particularly when the matrix is doubly stochastic \cite{nedic2018network,pu2021sharp,koloskova2021improved,alghunaim2022unified}. However, various alternative analytical strategies have emerged to characterize the performance of decentralized algorithms over more sophisticated network structures (directed graphs). Examples include (stochastic) gradient-push \cite{nedic2014distributed,nedic2016stochastic,assran2019stochastic,spiridonoff2020robust}, SONATA \cite{tian2020achieving,scutari2019distributed,kungurtsev2023decentralized}, Push-DIGing/ADD-OPT \cite{nedic2017achieving,liang2023understanding}, and Push-Pull/AB \cite{pu2020push,xin2018linear,xin2019distributed,xin2020general}.

In this paper, we propose a general analytical framework based on a set of infinite sums to study the convergence behavior of the Push-Pull algorithm under stochastic gradients (Stochastic Push-Pull). This novel approach enables a comprehensive convergence analysis for the algorithm, effectively covering various state-of-the-art results for decentralized stochastic gradient tracking (DSGT) methods beyond the limitations imposed by the spectral gap analysis. Additionally, we provide a sufficient condition under which the algorithm achieves linear speedup, matching the scalability of centralized stochastic gradient methods. Finally, we conduct numerical experiments that validate the theoretical findings.

\subsection{Related Works} 
Decentralized Stochastic Gradient Descent (DSGD) and its variants \cite{lian2017can,pu2021sharp,ying2021exponential,koloskova2019decentralized} have attracted considerable attention for large-scale distributed training. Although these methods are flexible and generally effective, they often encounter difficulties in handling data heterogeneity \cite{koloskova2020unified}. Such a challenge has motivated the development of more advanced algorithms, such as EXTRA \cite{shi2015extra}, Exact-Diffusion/$\text{D}^2$ \cite{li2019decentralized,tang2018d,huang2022improving,yuan2023removing}, and gradient tracking based methods \cite{pu2021distributed,koloskova2021improved}.

Most decentralized stochastic gradient algorithms rely on specific constraints on the weight matrices of the communication graphs (e.g., double stochasticity), which can be challenging or even impossible to satisfy for general directed graphs \cite{gharesifard2010does}. By contrast, methods such as subgradient-push (SGP) \cite{nedic2014distributed,nedic2016stochastic,assran2019stochastic}, SONATA \cite{tian2020achieving,scutari2019distributed} and Push-DIGing/ADD-OPT \cite{nedic2017achieving,xi2017add,liang2023understanding} utilize column-stochastic and/or row-stochastic weight matrices, enabling consensus optimization on strongly connected directed graphs. Several algorithms, including those in \cite{assran2019stochastic,spiridonoff2020robust,liang2023understanding}, have been shown to enjoy linear speedup when the number of iterations is sufficiently large.
In particular, the work in \cite{liang2023understanding} introduces a new metric that captures the influence of column-stochastic mixing matrices on decentralized stochastic gradient algorithms.

Another widely studied approach for optimization over general directed graphs is the Push-Pull/AB method \cite{xin2018linear, pu2020push,xin2020general}, originally designed for minimizing strongly convex objectives with deterministic gradients. The algorithm can operate using two separate subgraphs \cite{pu2020push}, each containing a (reversed) spanning tree that shares a common root. In recent years, an extensive family of Push-Pull based methods has emerged, addressing diverse scenarios. For example, the papers \cite{xin2019distributed,zhao2023asymptotic} investigates the performance of Push-Pull with stochastic gradients  under smooth strongly convex objectives, the work in \cite{zhu2024r} extends Push-Pull to the asynchronous setting, the papers \cite{saadatniaki2020decentralized,nguyen2023accelerated} introduce Push-Pull/AB variants to accommodate time-varying graphs, and \cite{song2022communication,liao2024robust} propose adaptations with compressed communication. Note that these methods often impose additional conditions, such as positive diagonal elements in the weight matrices. 

Despite the success of the Push-Pull method, its convergence behavior under stochastic gradients are less explored. Most existing works have not established comparable convergence guarantee for Stochastic Push-Pull against centralized SGD or push-sum based algorithms \cite{xin2019distributed,zhao2023asymptotic,nguyen2023accelerated,zhu2024r}. For example, while the work in \cite{zhao2023asymptotic} first demonstrates the $\mathcal{O}(\frac{1}{t})$ convergence rate of Push-Pull/AB with stochastic gradients under smooth strongly convex objectives, the dependency of the rate on the problem parameters remains unclear.
Recent works \cite{you2024b, you2025distributed} demonstrate the linear-speedup property of Stochastic Push-Pull with explicit transient iteration complexities, but the results heavily depend on global knowledge of the communication structure to construct appropriate spanning trees, which can be impractical for certain applications.

\subsection{Main Contribution}

In this work, we consider the Stochastic Push-Pull (S-PP) algorithm for minimizing smooth nonconvex objectives over general directed graphs. We begin by presenting an extended analysis of the underlying assumptions, particularly those related to the weight matrices, and then offer a detailed and novel convergence study. 

The main contribution of this paper is summarized as follows:
\begin{itemize}
    \item We provide a detailed discussion of the underlying assumptions for the Stochastic Push-Pull (S-PP) algorithm, with particular emphasis on the network structure and weight matrices, which generalizes from existing works. 
        \item Our analysis introduces a general analytical framework applicable to a broad class of distributed stochastic gradient algorithms. Unlike traditional approaches relying on the spectral gap, our framework leverages a series of finite sums to precisely characterize algorithmic convergence. This approach is particularly effective for analyzing the performance on directed and unbalanced graphs such as Erdos-Renyi graphs \cite{nedic2018network} and Multi-Sub-Ring graphs \cite{you2025distributed}.
    \item We show that, under smooth nonconvex objectives, S-PP generally achieves $\mathcal{O}(\frac{1}{\sqrt{T}})$ convergence rate with explicit coefficients. Notably, it is demonstrated for the first time that decentralized stochastic gradient tracking (DSGT) method achieves linear speedup with doubly stochastic (but not necessarily symmetric) mixing matrix which is compatible with balanced directed graphs (see Remark \ref{remark:case 1}). 
    For undirected graphs where symmetric mixing matrix is applicable, the convergence results of S-PP recover and/or enhance the state-of-the-art results for DSGT \cite{alghunaim2022unified,koloskova2021improved} (see Remark \ref{remark:case 23}). 
        \item  We establish a sufficient condition under which the S-PP algorithm achieves linear speedup, matching the scalability of centralized stochastic gradient methods--a result not previously reported. Such a condition is applicable to several existing stochastic gradient tracking methods that use different strategies for constructing the communication graphs and weight matrices, including those in \cite{alghunaim2022unified,koloskova2021improved,you2024b,you2025distributed}. 
     \item Numerical experiments demonstrate that S-PP is highly effective and robust across various decentralized optimization scenarios.
\end{itemize}

\subsection{Notation}
In this paper, all vectors are column vectors by default. Each node $i$ maintains a local copy $x_i \in \reals^p$ of the decision variable, an auxiliary variable $y_i\in \reals^p$. We denote their values at iteration $t$ by $x_i^{(t)}$ and $y_i^{(t)}$, respectively, and define 
\[
\begin{aligned}
    \bX & = \brk{x_1, x_2, \cdots,x_n}^\T \in \reals^{n\times p}, \\
    \bY & = \brk{y_1, y_2, \cdots,y_n}^\T \in \reals^{n\times p}.
\end{aligned} 
\]
We also define the aggregate gradient as 
\[
\begin{aligned}
    \nabla F(\bX) :=  \brk{  \nabla f_1 (x_1),\nabla f_2 (x_2),  \cdots, \nabla f_n(x_n) }^\T \in \reals^{n\times p}.
\end{aligned}
\] 
Each agent $i$ is able to query a stochastic gradient $g_i:= g_i(x_i,\xi_i)$ given $x_i\in \reals^p$, where $\xi_i$ represents some independent random variable. Let 
\[
\begin{aligned}
    \bxi &:= \brk{\xi_1, \xi_2, \cdots, \xi_n}^\T\\
    \bG(\bX,\bxi) & := \brk{g_1(x_1,\xi_1),  \cdots, g_n(x_n,\xi_n)}^\T \in \reals^{n\times p}.
\end{aligned}
\]

Let $\mone$ denote the column vector of all ones.
For convenience, we use the notations $\nabla \bF^{(t)}:=\nabla \bF(\bX^{(t)})$ and $\bG^{(t)}:= \bG(\bX^{(t)},\bxi^{(t)})$.
We denote by $\langle a, b\rangle$ the inner product of two vectors $a,b\in \reals^p$. The symbols $\norm{\cdot}_2$ and $\norm{\cdot}_F$ refer to the spectral norm and the Frobenius norm of a matrix respectively; for vectors, these norms reduce to the Euclidean norm $\norm{\cdot}$. Throughout the paper, we adopt the convention that any matrix raised to the zeroth power is the identity matrix of the corresponding dimension, and any summation whose upper limit is smaller than its lower limit equals zero. Furthermore, we use $a\lesssim b$ to indicate that $a\le Cb$ for some numerical constant $C$.

We consider a network of $n$ workers that cooperate to solve Problem \eqref{eq:obj} over a directed graph $\cG = \prt{\cN,\cE}$. Here, $\cN = \crk{1,2,\cdots,n}$ is the set of nodes, and $\cE\subseteq \cN\times\cN$ is the set of directed edges. An edge  $(j,i)\in \cE$ means that node $j$ can directly send information to node $i$. 

Let $\bW = [w_{ij}]\in \reals^{n\times n}$ be a nonnegative matrix. Define the induced graph $\cG_{\bW}$ by letting $(i,j)\in \cE_{\bW}$ if and only if $w_{ji} > 0$. For each node $i$, we define its in-neighborhood and out-neighborhood by $\cN_{\bW,i}^{\text{in}}:= \crk{j| w_{ij}>0, j\in \cN}$ and $\cN_{\bW,i}^{\text{out}}:= \crk{j| w_{ji}>0, j\in \cN}$, respectively. 
Furthermore, let $\cR_{\bW}$ be the set of roots for all possible spanning trees in the graph $\cG_{\bW}$.

\subsection{Organization}
The remaining parts of this paper are organized as follows. In Section \ref{sec:setup}, we introduce the standing assumptions and the main algorithm. In Section \ref{sec:analysis}, we conduct the convergence analysis. Then, we present the convergence results in Section \ref{sec:convergence}. Section \ref{sec:experiment} showcases the numerical experiments. Finally, we conclude the paper in Section \ref{sec:conclusion}.

\section{Setup}
\label{sec:setup}

In this section, we first introduce the standing assumptions in Subsection \ref{subsec:assumption} with some necessary discussions. Then, we present the Stochastic Push-Pull (S-PP) algorithm in Subsection \ref{subsec:algo}.

\subsection{Assumptions}
\label{subsec:assumption}

We introduce the fundamental assumptions for the considered S-PP algorithm which also underpin the follow-up analysis. First, in Assumption \ref{a.graph1}, we specify the structures of the communication graphs, a requirement that remains consistent with the standard practice in the existing literature \cite{pu2020push,zhao2023asymptotic,zhu2024r}. Second, in Assumption \ref{a.graph}, we introduce the conditions on the weight matrices. Unlike previous approaches that demand strictly positive diagonal entries \cite{pu2020push,zhao2023asymptotic,zhu2024r}, we relax the requirement to an exponentially decaying property (Definition \ref{def:exp}), thereby broadening the class of admissible weighting schemes and enhancing the flexibility of network designs. Assumption \ref{a.var} and \ref{a.smooth} are standard in the analysis of stochastic gradient algorithms \cite{lian2017can,ying2021exponential,you2024b}.

We first specify the conditions for the two (possibly identical) graphs used in S-PP, following the paper \cite{pu2020push}. The two graphs $\cG_\bR$ and $\cG_{\bC}$ correspond to two matrices $\bR$ and $\bC$, respectively, which will be discussed later.
\begin{assumption}[Communication Graphs]\label{a.graph1}
    The graphs $\cG_\bR$ and $\cG_{\bC^\T}$\footnote{The graph $\cG_{\bC^\T}$ is essentially $\cG_{\bC}$ with all edges reversing directions.} each contain at least one spanning tree, and there exists at least one pair of spanning trees (one in $\cG_\bR$ and one in $\cG_{\bC^\T}$) that share a common root. Equivalently, there hold $\cR := \cR_{\bR} \cap \cR_{\bC^\T} \neq \emptyset$ and $r := |\cR| \ge 1$.
\end{assumption}

We now introduce the definitions for characterizing the condition on the weight matrices.
\begin{definition}[Root Eigenvector]
    \label{def:root}
    Let $\bA \in \mathbb{R}^{n\times n}$ be a nonnegative, row-stochastic matrix. A nonnegative, unit left eigenvector $\bpi_\bA$ associated with the eigenvalue~$1$ (i.e., $\bpi_\bA^\T \bA = \bpi_\bA^\T$, $\bpi_\bA \ge 0$, and $\|\bpi_\bA\|_1 = 1$) is called a \emph{root eigenvector} of $\bA$ if the $i$-th element $\bpi_{\bA}(i) = 0$ for all $i \notin \cR_{\bA}$.
\end{definition}
    
\begin{definition}[Exponential Decay]
    \label{def:exp}
    A nonnegative, row-stochastic matrix $\bA \in \mathbb{R}^{n\times n}$ is said to be \emph{exponentially decaying} if there exist a constant $\alpha \in (0,1)$, 
    an integer $m > 0$, and a root eigenvector $\bpi_\bA$ such that, for all $t \ge m$,
    \[
        \|\bA^{t} - \mone\bpi_\bA^\T\|_2 \le \alpha^{t}.
    \]
\end{definition}
\begin{remark}
    Observe that for a doubly stochastic matrix $\bW$ satisfying $\norm{\bW - \mone\mone^\T/n}_2 < 1$, the vector $\mone/n$ is both a root eigenvector and a right eigenvector corresponding to the eigenvalue $1$. Consequently, $\bW$ exhibits an exponentially decaying property with parameters $m=1$ and $\alpha = \norm{\bW - \mone\mone^\T/n}_2$. For simplicity, we omit explicitly stating this condition when later referring to a doubly stochastic matrix $\bW$ satisfying Assumption \ref{a.graph}.
\end{remark}

\begin{assumption}[Weight Matrices]
    \label{a.graph}
    The matrix $\bR$ is row-stochastic and the matrix $\bC$ is column-stochastic, i.e., $\bR \mone = \mone$ and $\bC^\T \mone = \mone$. In addition, $\bR$ and $\bC^\T$ are exponentially decaying.
\end{assumption}

From the assumptions above, we obtain the following critical result related to the information flow among the nodes.
\begin{lemma}
    \label{l:eigen}
    Under Assumptions \ref{a.graph1} and \ref{a.graph}, the weight matrix $\bR$ has a unique root eigenvector $\bpi_\bR^\T$, and $\bC^\T$ has a unique root eigenvector $\bpi_\bC^\T$. Moreover, $\bpi_\bR^\T \bpi_\bC > 0$.
\end{lemma}
    
\begin{proof}
    See Appendix \ref{pf:l:eigen}.
\end{proof}
    Assumptions \ref{a.graph1} and \ref{a.graph} establish a minimal requirement to ensure unobstructed information flow among the nodes. In contrast, many existing algorithms (e.g., \cite{lian2017can,xin2018linear,ying2021exponential,li2019decentralized,nedic2018network}) require the communication graph to be strongly connected, limiting the flexibility in the architecture design. Moreover, earlier works related to Push-Pull (e.g., \cite{pu2020push,zhao2023asymptotic,zhu2024r}) often require the weight matrices $\bR$ and $\bC$ to have strictly positive diagonal entries rather than the exponential decay property introduced in Assumption \ref{a.graph}. We formally show in Lemma \ref{lem:diagonal} below that exponential decay is a weaker condition. Indeed, it allows more flexible choices for $\bR$ and $\bC$. Examples of mixing matrices satisfying Assumption \ref{a.graph} but may not have strictly positive diagonal entries include, for instance, arbitrary doubly stochastic matrices, 
    $(0,1)$-mixing matrices induced by $B$-ary tree families \cite{you2024b} or arbitrary spanning trees \cite{you2025distributed}.

\begin{lemma}
    \label{lem:diagonal}
    Given a non-negative, row-stochastic matrix $\bA$ with positive diagonal elements. Assume its induced graph $\cG_{\bA}$ contains at least one spanning tree. Then, the matrix $\bA$ is exponential decaying with a root eigenvector denoted by $\bpi^\T$. More specifically, there exists an positive integer $m_{\bA}$ related only  to $\bA$, such that, for any $m\ge m_{\bA}$, it holds that
    \[
    \norm{\bA^m - \mone\bpi^\T}_2 \le (1-p)^{m},
    \]
    where $p = \frac{1}{2}\prt{1- \lambda}\in (0,1)$ and $\lambda$ is the second largest eigenvalue of $\bA$ w.r.t. modulus.
\end{lemma}
\begin{proof}
    See Appendix \ref{pf:lem:diagonal}.
\end{proof}

Regarding the stochastic gradients, we consider the following standard assumption.
\begin{assumption}
    \label{a.var}
    Each node $i \in \cN$ is able to obtain noisy gradient $g(x,\xi_i)$ given $x\in \reals^p$, where each random vector $\xi_i \sim \cD_i$ is independent across $i\in \cN$. In addition, for some $\sigma^2 > 0$, we have
    \[
    \begin{aligned}
    & \expect_{\xi_i\sim\cD_i} \crk{ g_i(x,\xi_i)| x}  = \nabla f_i(x) , \\
    & \expect_{\xi_i \sim \cD_i} \crk{ \norm{g_i(x,\xi_i) - \nabla f_i(x)}^2 |x}  \le \sigma^2.
    \end{aligned}
    \]
\end{assumption}
Stochastic gradients appear in various areas including online distributed learning, reinforcement learning, generative modeling, and parameter estimation; see, e.g., \cite{you2024b,lian2017can,pu2021distributed,ying2021exponential}. Moreover, Assumption \ref{a.var} generally holds for the empirical risk minimization (ERM) problems if the stochastic gradients are queried through uniformly sampling with replacement \cite{huang2024distributed}.

Assumption \ref{a.smooth} below is common that requires the individual ojective functions $f_i$ to be smooth and lower bounded.
\begin{assumption}
    \label{a.smooth}
        Each $f_i(x) : \reals^p \rightarrow \reals$ is lower bounded with $L$-Lipschitz continuous gradients, i.e. for any $x, x' \in \reals^p$,
        \[
        \norm{ \nabla f_i(x) - \nabla f_i(x')} \le L\norm{x - x'}.
        \]
        We denote $f^* := \min_{x\in \reals^p} f(x)$.
\end{assumption}

\subsection{Algorithm}
\label{subsec:algo}
In this subsection, we introduce the Stochastic Push-Pull (S-PP) algorithm which is formally outlined in Algorithm \ref{alg:pp}.\footnote{The Push-Pull/AB type methods with stochastic gradients have been proposed and studied under different names (e.g., \emph{SAB} in \cite{xin2019distributed}) with more restricted assumptions. We use \emph{Stochastic Push-Pull} to emphasize the connection to the original Push-Pull method \cite{pu2020push} that first considers Assumption \ref{a.graph1}.} The method allows each agent to communicate with its neighbors over two communication graphs, $\cG_\bR$ (Pull Graph) and $\cG_{\bC}$ (Push Graph) satisfying Assumption \ref{a.graph1} and \ref{a.graph}, and solve Problem \eqref{eq:obj} in a distributed manner. 

\begin{algorithm}[htbp]
    \caption{Stochastic Push-Pull Algorithm (\textbf{S-PP})}
    \label{alg:pp}
    \begin{algorithmic}[1]
        \Require Communication graphs $\cG_\bR$ (Pull Graph) and $\cG_{\bC}$ (Push Graph). Each agent $i$ initializes with any arbitrary but identical $x_i^{(0)} = x^{(0)} \in \reals^p$, stepsize $\gamma$, initial gradient tracker $y_i^{(0)} = g_i(x_i^{(0)},\xi_i^{(0)})$ by drawing a random sample $\xi_i^{(0)}$. 
        \For{Iteration $t = 0, 1, 2,\ldots, T-1$}
        \For{Agent $i$ \textbf{in parallel}}
        \State Pull $x_j^{(t)} - \gamma y_{j}^{(t)}$ from each $j \in \cN_{\bR,i}^{\text{in}}$.
        \State Push $y_{j}^{(t)}$ to each $j\in \cN_{\bC,i}^{\text{out}}$.
        \State Independently draw a random sample $\xi_i^{(t+1)}$ .
        \State \textbf{Update} parameters through
        \[
        \begin{aligned}
            x_i^{(t+1)} & = \sum_{j \in \cN_{\bR,i}^{\text{in}} } r_{ij}\prt{x_j^{(t)} - \gamma y_{j}^{(t)}}, \\
            y_i^{(t+1)} & = \sum_{j\in\cN_{\bC,i}^{\text{in}} } c_{ij} \prt{y_j^{(t)} + g_i^{(t+1)} - g_i^{(t)} }.
        \end{aligned}
        \]
        \EndFor
        \EndFor
    \end{algorithmic}
\end{algorithm}
S-PP consists of two major steps: a pull step and a push step. In the pull step, each agent $i$ pulls the model information through the graph $\cG_\bR$ from its in-neighbors $\cN_{\bR,i}^{\text{in}}$ and updates its local variable $x_i$ accordingly. In the push step, agent $i$ pushes the information regarding the stochastic gradients to its out-neighbors $\cN_{\bC,i}^{\text{out}}$ in the graph $\cG_\bC$ and updates its auxiliary variable $y_i$. The algorithm iterates over the pull and push steps until a stopping criterion is met. As will be shown in Section \ref{sec:convergence}, S-PP is guaranteed to converge to a stationary point of the objective function in expectation under certain conditions.

For the ease of presentation and analysis, we often consider the following compact form of Algorithm \ref{alg:pp}.
\begin{equation}
    \label{eq:pp}
    \begin{aligned}
        \bX^{(t+1)} & = \bR\bX^{(t)} - \gamma\bR\bY^{(t)}, \\
        \bY^{(t+1)} & = \bC \prt{\bY^{(t)} + \bG^{(t+1)} - \bG^{(t)}}.
    \end{aligned}
\end{equation}

\section{Preliminary Results}
\label{sec:analysis}
In this section, we present several preliminary results focusing on three aspects: characterizing the influence of directed graphs, introducing useful analytical tools, and convergence analysis. 

We define the following matrices related to the weight matrices $\bR$ and $\bC$ satisfying Assumptions \ref{a.graph1} and \ref{a.graph}:
\[
\bPi_{\bR} := \bI - \mone\bpi_{\bR}^\T,  \bPi_{\bC} := \bI - \bpi_{\bC}\mone^\T .
\]
For simplicity, we define $\tilde{\bR} := \bPi_{\bR} \bR$, $\tilde{\bC}:= \bPi_{\bC}\bC$ and $\pi = \bpi_\bR^\T\bpi_\bC$. Furthermore, let $\tilde{\bR}^{k} := \bPi_{\bR}\bR^{k}$  and $\tbC^{k}:= \bPi_{\bC}\bC^{k}$  for any integer $k\ge 0$.

 \subsection{Influence of the Directed Graphs}
Unlike prior analyses \cite{pu2020push,zhao2023asymptotic,zhu2024r} which rely on specific matrix norms tailored to matrices $\bR$ and $\bC$, we instead consider a set of infinite sums that capture how efficiently information flows converge at the common root nodes in $\cR$ and propagate through $\cG_{\bR}$ and $\cG_{\bC}$. As will be shown in Section \ref{sec:convergence}, these infinite-sum characterizations also help recover known results for existing decentralized stochastic gradient tracking (DSGT) methods on undirected graphs \cite{alghunaim2022unified, koloskova2021improved}. 

Below, we first define the sums $M_1$ and $M_2$ which play a crucial role in studying the linear speedup property for the Stochastic Push-Pull algorithm:
\[
\begin{aligned}
    M_1 & := \norm{\bpi_\bR^\T \bC}_2^2 + \sum_{t=1}^{\infty}\norm{\bpi_{\bR}^\T\prt{\bC^{t+1} - \bC^{t}}}_2^2, \ M_2 := \sum_{t=1}^{\infty} t\norm{\bpi_{\bR}^\T\prt{\bC^{t+1} - \bC^{t}}}_2.
\end{aligned}
\]
We then define the sums $N_1$ through $N_8$, which provide a more refined characterization of the convergence rate in our analytical framework.
\[
\begin{aligned}
    N_1 & := \sum_{t=1}^{\infty}\norm{\bpi_\bR^\T  \tbC^{t}}_2 ,\ N_2 := \sum_{t=0}^{\infty}\norm{\bpi_\bR^\T \tbC^{t}}_2^2,\ N_3  := \sum_{t=1}^{\infty}\norm{\tbR^t\bpi_\bC}_2,\ N_4 := \sum_{t=1}^{\infty}\norm{\tbR^t\bpi_\bC}_2^2,\\
    N_5 &:= \sum_{t=1}^{\infty} \norm{\sum_{k=1}^{t-1}\tbR^k \tbC^{t-k}}_2, \ N_6 := \sum_{t=1}^{\infty} \norm{\sum_{k=1}^{t} \tbR^k \tbC^{t-k}}_2^2,\\
    N_7 & :=  \sum_{t=1}^{\infty} \norm{\sum_{k=1}^{t}\tbR^k \tbC^{t-k+1} - \sum_{k=1}^{t-1}\tbR^k \tbC^{t-k}}_2^2, \ N_8 := \sum_{t=1}^{\infty} \norm{\sum_{k=1}^{t}\tbR^k \tbC^{t-k+1} - \sum_{k=1}^{t-1}\tbR^k \tbC^{t-k}}_2.
\end{aligned}
\]

Lemma \ref{thm:series} below shows that all the sums above are finite.
\begin{lemma}
    \label{thm:series}
    Suppose that $\bR$ and $\bC$ satisfy Assumptions \ref{a.graph1} and \ref{a.graph}. Then, for all $i=1,2$ and $j=1,\cdots,8$, the sums $M_i$ and $N_j$ are finite, i.e.,
    \[
    M_i < \infty, \text{ and } N_j < \infty.
    \]
\end{lemma}
\begin{proof}
    The general idea is to split each infinite series into two parts: one that converges to a finite value and the other one with a finite number of terms. As a first step, we show that the series $M_1$ converges (i.e., bounded). 
    
    By Assumption \ref{a.graph}, $\bR$ and $\bC$ are exponentially decaying, i.e., there exists an integer $m_{\bR,\bC}:=\max\crk{m_\bR, m_\bC}$ such that for any $t > m_{\bR,\bC}$,  
\[
\norm{\tbR^t }_2 \le (1-p)^{t}, \norm{\tbC^t }_2 \le (1-p)^{t},
\]
where $p < 1$.
Then, we denote
\[
\begin{aligned}
    N := \max & \left\{ \norm{\bpi_\bR}_2,\norm{\bpi_\bC}_2, \norm{\tbR^t }_2, \norm{\tbC^t }_2, \prt{1-p}^{m_{\bR,\bC}+1},\forall t \le  m_{\bR,\bC} \right\}. 
\end{aligned}
\]
For $M_1$, notice that 
\[
\bC^{t+1} - \bC^{t} = \prt{\bC^{t+1} - \bpi_\bC\mone^\T} - \prt{\bC^{t} - \bpi_\bC\mone^\T}.
\] 
From the inequality that $\norm{a+b}_2^2 \le 2\norm{a}_2^2 + 2\norm{b}_2^2$, we have
\[
\begin{aligned}
    & \norm{\bpi_{\bR}^\T\prt{\bC^{t+1} - \bC^{t}}}_2^2  \le 2 \norm{ \bpi_\bR\tbC^{t+1}  }_2^2 + 2 \norm{ \bpi_\bR\tbC^{t}  }_2^2 \le  2 \norm{\bpi_\bR}_2^2\norm{ \tbC^{t+1}   }_2^2 + 2 \norm{\bpi_\bR}_2^2\norm{ \tbC^{t} }_2^2.
\end{aligned}
\]
Then, it holds that 
\[
\norm{\bpi_{\bR}^\T\prt{\bC^{t+1} - \bC^{t}}}_2^2 \le \left\{    
\begin{aligned} 
    & 4N^4 & \text{ for } t \le m_{\bR,\bC}, \\
    & 4N^2\prt{1-p}^{2t} & \text{ for } t > m_{\bR,\bC}.
\end{aligned} 
\right.
\]
Summing over $t$ from $1$ to infinity, we have
\[
\begin{aligned}
    & \sum_{t=1}^{\infty} \norm{\bpi_{\bR}^\T\prt{\bC^{t+1} - \bC^{t}}}_2^2 \le 4 N^4 m_{\bR,\bC} + 4N^2 \sum_{t=1}^{\infty} \prt{1-p}^{2t}  \le  4 N^4 m_{\bR,\bC} + 4N^2\frac{(1-p)^2}{2p-p^2} < \infty.
\end{aligned} 
\]
Then, we conclude that
\[
M_1 \le \norm{\bpi_\bR^\T \bC}_2^2 + 4 N^4 m_{\bR,\bC} + 4N^2\frac{(1-p)^2}{2p-p^2} < \infty.
\]
Proofs for the remaining terms follow similar arguments by using the identity $\sum_{t=1}^{\infty} t(1-p)^t = (1-p)/p^2$. Hence, we omit those details.
\end{proof}

We now discuss the scenario $\bR = \bC = \bW$, where $\bW$ is nonnegative and doubly stochastic, satisfying Assumption \ref{a.graph} with $\lambda:=\norm{\bW - \mone\mone^\T/n}_2 < 1$. Also, assume Assumption \ref{a.graph1} holds for $\cG_{\bW}$. Such a setting is commonly considered in the literature for undirected graphs \cite{alghunaim2022unified,koloskova2021improved, pu2021distributed} and Lemma \ref{l:symmetric} applies in these cases.

\begin{lemma}
    \label{l:symmetric}
    For $\bR = \bC = \bW$, where $\bW$ is nonnegative, doubly stochastic and satisfying Assumption \ref{a.graph}. Additionally, assume Assumptions \ref{a.graph1} holds for $\cG_{\bW}$, and define $\lambda:=\norm{\bW - \mone\mone^\T/n}_2 < 1$, we have
    \[
    \begin{aligned}
        M_1 & = 1/n, M_2 = N_1 = N_2 = N_3 = N_4 = 0, N_5 \le \frac{\lambda^2}{(1-\lambda)^2},\\
        N_6 & \le \frac{\lambda^2 + \lambda^4}{(1-\lambda^2)^3}, N_7 \le \frac{8}{\prt{1-\lambda^2}^3}, N_8  \le  \frac{8}{(1 - \lambda)^2}, 
    \end{aligned}
    \]
    Furthermore, if $\bW$ is additionally symmetric, it holds exactly for
    \[
    \begin{aligned}
        N_5 & = \frac{\lambda^2}{(1-\lambda)^2}, N_6 = \frac{\lambda^2 + \lambda^4}{(1-\lambda^2)^3}, N_7  \le \frac{10}{\prt{1 + \min\crk{\lambda_n,0}}\prt{1-\lambda}}, N_8  \le  \frac{10}{\sqrt{1 + \min\crk{\lambda_n,0}}(1 - \lambda)}, 
    \end{aligned}
    \]
    where $\lambda_n \in (-1,1)$ is the smallest eigenvalue of $\bW$.
\end{lemma}
\begin{proof}
    See Appendix \ref{pf:l:symmetric}.
\end{proof}
It can be seen that simple characterizations for the constants are available in such a case.

\subsection{Supporting Lemmas}

In this subsection, we present several inequalities and lemmas that are instrumental for proving the main results.
\begin{lemma}
    \label{lem:sum_matrix}
    For an arbitrary set of $p$ matrices $\crk{\bA_i}_{i=1}^{p}$ with the same size, and given $\alpha_1,\cdots, \alpha_p$ satisfying $\alpha_i > 0$ and $\sum_{i=1}^{p}\alpha_i \le 1$, we have
    \[
    \norm{\sum_{i=1}^{p}\bA_i}_F^2 \le \sum_{i=1}^{p} \frac{1}{\alpha_i} \norm{\bA_i}_F^2.
    \]
\end{lemma}
\begin{proof}
    See \cite[Lemma E.1]{you2025distributed} for reference.
\end{proof}

\begin{lemma}
    \label{lem:matrix_norm}
    Let $\bA$, $\bB$ be two matrices whose sizes match. Then,
    \[
    \norm{\bA\bB}_F \le \norm{\bA}_2\norm{\bB}_F.
    \]
\end{lemma}
\begin{proof}
See \cite[Lemma A.6]{you2024b} for reference.
\end{proof}

\begin{lemma}
    \label{lem:martin}
    Let $A,B,C$ and $\alpha$ be positive constants and $T$ be a positive integer. Define 
    \[
    g(\gamma) = \frac{A}{\gamma (T+1)} + B\gamma + C\gamma^2.
    \]
    Then,
    \[
    \inf_{\gamma \in (0,\frac{1}{\alpha}]} g(\gamma) \le 2\prt{\frac{AB}{T+1}}^{\frac{1}{2}} + 2C^{\frac{1}{3}} \prt{\frac{A}{T+1}}^{\frac{2}{3}} + \frac{\alpha A}{T+1},
    \]
    where the upper bound can be achieved by choosing $\gamma = \min\crk{\prt{\frac{A}{B(T+1)}}^{\frac{1}{2}} , \prt{\frac{A}{C(T+1)}}^{\frac{1}{3}},\frac{1}{\alpha}}$.
\end{lemma}
\begin{proof}
    See \cite[Lemma 26]{koloskova2021improved} for reference.
\end{proof}

\begin{lemma}
    \label{lem:sum_help}
    Given non-negative sequences $\crk{a_i}_{i=0}^{\infty}$ and $\crk{b_i}_{i=0}^{\infty}$, suppose that $\sum_{i=0}^{\infty} a_i < \infty$. Then, for any integer $T>0$, we have
    \[
    \sum_{t=0}^{T} \sum_{j=0}^{t} a_{t-j} b_j \le \sum_{i=0}^{\infty} a_i \sum_{j=0}^{T} b_j.
    \]
\end{lemma}
\begin{proof}
    The inequality holds directly by changing the order of summation.
\end{proof}

\begin{lemma}
    \label{lem:series_help1}
    Given a positive series $\crk{a_n}_{n=1}^{\infty}$ and integers $T,k$ s.t. $T>k>0$, it holds that
    \[
    \sum_{t=k+1}^{T}\sum_{j=1}^k a_{t-j} \le \sum_{i=1}^{T-1} i a_i.
    \]
\end{lemma}
\begin{proof}
    See Appendix \ref{pf:lem:series_help1}.
\end{proof}

\subsection{Convergence Analysis}
We now present the key steps for analyzing the convergence of S-PP. We begin by defining several notations for simplicity:
\[
\begin{aligned}
    \hat{x}^{(t)} & := \bpi_\bR^\T\bX^{(t)},\Bar{\bX}^{(t)} := \mone \hat{x}^{(t)} , \Bar{\bY}^{(t)}:= \bpi_{\bC}\mone^\T\bY^{(t)} \\
    \theta^{(t)}_i & := g_i^{(t)} - \nabla f_i(x_i^{(t)}), \btheta^{(t)}  := \bG^{(t)} - \nabla \bF(\bX^{(t)}), \\
     \hat{\bX}^{(t)} &:= \bPi_\bR \bX^{(t)}, \Delta \bar{\bX}^{(t)}  : = \bar{\bX}^{(t+1)} -  \bar{\bX}^{(t)},\nabla \bbF^{(t)}  := \nabla \bF(\bX^{(t)}) - \nabla \bF(\bar{\bX}^{(t)}).
\end{aligned}
\]
Here, $\hat{x}^{(t)}$ represents the desired output of Algorithm \ref{alg:pp} at iteration $t$.

Next, we express the S-PP algorithm, as outlined in Algorithm \ref{alg:pp}, in a succinct matrix form. From Equation \eqref{eq:pp}, the update rules for $\bX^{(t)}$ and $\bY^{(t)}$ are:
\begin{equation}
    \label{eq:pp_matrix}
    \begin{aligned}
        \prt{\begin{matrix}
            \bX^{(t+1)} \\
            \bY^{(t+1)}
        \end{matrix}} = &  \prt{\begin{matrix}
            \bR  & - \gamma\bR   \\
            \mathbf{0} & \bC 
        \end{matrix}} \prt{\begin{matrix}
            \bX^{(t)} \\
            \bY^{(t)}
        \end{matrix}}  + \prt{\begin{matrix}
            \mathbf{0}  & \mathbf{0}   \\
            \mathbf{0} & \bC 
        \end{matrix}} \prt{\begin{matrix}
            \mathbf{0}\\
            \bG^{(t+1)} - \bG^{(t)}
        \end{matrix}}.
    \end{aligned}
\end{equation}
It follows that
\[
\begin{aligned}
    & \prt{\begin{matrix}
         \bX^{(t)}\\
        \bY^{(t)}  
    \end{matrix}}
     =  \prt{\begin{matrix}
         \bR & -\gamma \bR\\
        \mathbf{0} & \bC  
    \end{matrix}}^t\prt{\begin{matrix}
         \bX^{(0)}\\
        \bY^{(0)}  
    \end{matrix}}  + 
    \sum_{j=0}^{t-1}
    \prt{\begin{matrix}
         \bR & -\gamma \bR\\
        \mathbf{0} & \bC  
    \end{matrix}} ^{t-j-1}\prt{\begin{matrix}
        \mathbf{0}  & \mathbf{0}   \\
        \mathbf{0} & \bC 
    \end{matrix}}\prt{\begin{matrix}
         \mathbf{0}\\
        \bG^{(j+1)} -\bG^{(j)} 
    \end{matrix}}.
\end{aligned}
\]
Note that for any integer $j>0$, we have the following matrix identity:
\[
\left( \begin{array}{cc}
         \bR & -\gamma \bR\\
        \mathbf{0} & \bC  
    \end{array}\right)^j = \left( \begin{array}{cc}
         \bR^j & -\gamma \sum_{k=1}^{j}\bR^k\bC^{j-k}\\
        \mathbf{0} & \bC^j  
    \end{array}\right).
\]
Therefore, in light of the initial condition $\bY^{(0)} = \bG^{(0)}$, we can express $\bX^{(t)}$ and $\bY^{(t)}$ as follows.
\begin{align}
    \bX^{(t)} =  &\bR^{t}\bX^{(0)} - \gamma \sum_{j=0}^{t-2} \sum_{k=1}^{t-j-1}\bR^{k}\bC^{t-j-k} \brk{\bG^{(j+1)} - \bG^{(j)}} -\gamma \sum_{k=1}^{t}\bR^k\bC^{t-k} \bG^{(0)}, \label{eq:pp1} \\
    \bY^{(t)}  = & \sum_{j=0}^{t-1} \bC^{t-j} \brk{\bG^{(j+1)} - \bG^{(j)}} + \bC^{t} \bG^{(0)}. \label{eq:pp2} 
\end{align}
After multiplying both sides of Equation \ref{eq:pp2} by $\bPi_{\bC}$, we obtain the following frequently applied transformation.
\begin{equation}
    \label{eq:pp2_PiC}
    \begin{aligned}
        \bPi_{\bC}\bY^{(t)} & = \sum_{j=0}^{t-1} \tbC^{t-j} \brk{\bG^{(j+1)} - \bG^{(j)}} + \tbC^{t} \bG^{(0)}  = \sum_{j=1}^{t-1}\prt{\tbC^{t-j+1} - \tbC^{t-j}} \bG^{(j)} + \tbC \bG^{(t)}.
    \end{aligned}
\end{equation}

To support the convergence analysis, we introduce several lemmas that play a crucial role in deriving the main results. Let $\cF_t$ be the $\sigma$-algebra generated by $\crk{\xi^{(j)}_i,i\in \brk{n}, 0\le j\le t-1}$, and define $\expect\brk{\cdot|\cF_t}$ as the conditional expectation given $\cF_t$. Lemma \ref{l:variance} provides an estimate for the variance of the gradient estimator $G(\bX^{(t)}, \bxi^{(t)})$. 

\begin{lemma}
    \label{l:variance}
    Under Assumption \ref{a.var}, given a vector $a\in \reals^{n}$, we have for all $t\ge 0$ that
    \[
    \begin{aligned}
        \expect\brk{\norm{a^\T \btheta^{(t)}}^2} & \le \norm{a}^2\sigma^2.
    \end{aligned}
    \]
\end{lemma}
\begin{proof}
    See Appendix \ref{pf:l:variance}.
\end{proof}
Lemma \ref{l:co-var} provides the covariance between $\bG^{(t-i)}$ and $\bG^{(t-j)}$ for delays $i$ and $j$.
\begin{lemma}
    \label{l:co-var}
    Suppose Assumption \ref{a.var} holds. Then, for any vector $a,b\in \reals^n$ and integer $i\ne j$, we have
    \[
    \expect\left\langle a^\T \btheta^{(t-i)}, b^\T \btheta^{(t-j)} \right\rangle = 0.
    \]
\end{lemma}

The following lemmas delineate the critical elements for constraining the average expected norms of the objective function formulated in Equation \eqref{eq:obj}, specifically, $\frac{1}{T+1}\sum_{t=0}^{T} \expect\norm{\nabla f (\hat{x}^{(t)})}^2$. In particular, Lemma \ref{l:X_diff} and Lemma \ref{l:PiX} provide the upper bounds for $\sum_{t=0}^{T}\expect\norm{\Delta\bar{\bX}^{(t)} }_F^2$ and $\sum_{t=0}^{T}\expect\norm{\hat{\bX}^{(t)}}_F^2$ respectively.
\begin{lemma}
    \label{l:X_diff}
    Suppose Assumption \ref{a.graph1}, \ref{a.graph}, \ref{a.var} and \ref{a.smooth} hold. For $\gamma \le \frac{1}{6\sqrt{n} N_1 L} $, we have the following inequality:
    \[
    \begin{aligned}
    & \sum_{t=0}^{T} \expect\norm{\Delta\bar{\bX}^{(t)}  }_F^2  \le 4\gamma^2 n M_1 \sigma^2 (T+1)  + 96 \gamma^2 n \max\crk{N_1^2, n \pi^2} L^2 \sum_{t=0}^{T} \expect\norm{\hat{\bX}^{(t)}}_F^2 \\
    &  + 24 \gamma^2  n^3\pi^2 \sum_{t=0}^{T} \expect\norm{\nabla f(\hat{x}^{(t)})}^2  + 12 \gamma^2 n N_2
      \norm{\nabla \bF^{(0)}}_F^2.
    \end{aligned}
    \]
\end{lemma}
\begin{proof}
    See Appendix \ref{pf:l:X_diff}.
\end{proof}
\begin{lemma}
    \label{l:PiX}
    Suppose Assumption \ref{a.graph1}, \ref{a.graph}, \ref{a.var} and \ref{a.smooth} hold. Then, for $\gamma \le \frac{1}{20\sqrt{C_1} L}$, where 
    \[
    C_1 =  \max\crk{nN_1^2, nN_3^2, N_8^2, \sqrt{n}N_1N_5, n\pi N_5},
    \]
    we have
    \[
\begin{aligned}
    & \sum_{t=0}^{T}\expect\norm{\hat{\bX}^{(t)} }_F^2 \le 864\gamma^2 \max\crk{N_4,N_7} n \sigma^2 \prt{T+1}+ 2400\gamma^4 nM_1N_5^2L^2 \sigma^2 (T+1) \\
    & +  210 \gamma^2 \max\crk{ N_2 N_5/\pi , N_6 } \norm{\nabla F^{(0)}}_F^2  + 240 \gamma^2 n \max\crk{n N_3^2, n\pi  N_5 } \sum_{t=0}^{T} \expect\norm{\nabla f(\hat{x}^{(t)})}^2.
\end{aligned}
    \]
\end{lemma}
\begin{proof}
    See Appendix \ref{pf:l:PiX}.
\end{proof}

To derive the main convergence result, we integrate the findings from Lemma \ref{thm:series}, Lemma \ref{l:X_diff} and Lemma \ref{l:PiX}, as outlined in Lemma \ref{l:descent} below.
\begin{lemma}
    \label{l:descent}
    Suppose Assumption \ref{a.graph1}, \ref{a.graph}, \ref{a.var} and \ref{a.smooth} hold. Define the following constants:
    \[
    \begin{aligned}
        P_1 = & \max\crk{nN_1^2, nN_3^2,N_8^2, \sqrt{n}  N_1 N_5, n\pi N_5} , P_2 =  \max\crk{n^2\pi^2, nN_1^2,  nN_1^2 \tilde{M}_2, n^2\pi^2 \tilde{M}_2^2} ,\\
        P_3 = & \max\crk{  n^2 \pi^2, \frac{N_1^4}{\pi^2},\frac{N_3^4}{\pi^2}},  P_4 =  \sqrt{\max\crk{P_2,P_3}}N_5, Q = \max\crk{M_1,\tilde{M}_2, M_1 \tilde{M}_2}.
    \end{aligned}
    \]
    Then, for $\gamma \le \frac{1}{500\sqrt{\max\crk{P_1,P_2,P_3,P_4}}L}$, we have
    \[
    \begin{aligned}
        & \frac{1}{T+1}\sum_{t=0}^{T} \expect\norm{\nabla f(\hat{x}^{(t)})}_2^2 \le \frac{10\Delta_f}{\gamma n \pi \prt{T+1}} + \frac{100Q\sigma^2 L}{n \pi } \gamma + \frac{ 80000\sqrt{\max\crk{P_2,P_3}} \max\crk{N_4, N_7} \sigma^2 L^2  }{n \pi }\gamma^2 \\
        & + \frac{200000\sqrt{\max\crk{P_2,P_3}} M_1 N_5^2\sigma^2 L^4 }{n \pi }\gamma^4  + \frac{30 P_5}{n^2\pi\prt{T+1}}F_0,
    \end{aligned}
    \]
    where $P_5 := \max\crk{N_2/\pi,N_2/(n\pi^2) , N_6/N_5}$, $F_0 = \norm{\nabla\bF^{(0)}}_F^2/n$ and $\tilde{M}_2=0$ if $\cG_{\bR}$ and $\cG_{\bC^\T}$ are the spanning trees considered in \cite{you2024b,you2025distributed} and $\tilde{M}_2=M_2$ otherwise.
\end{lemma}
\begin{proof}
    See Appendix \ref{pf:l:descent}.
\end{proof}

\section{Convergence Results}
\label{sec:convergence}
In this section, we present the main convergence results for the Stochastic Push-Pull algorithm, detailing its convergence behavior. Notably, these results cover the case when $\bR = \bC = \bW$, where $\bW$ is symmetric and doubly stochastic, corresponding to the distributed stochastic gradient tracking (DSGT) method studied in \cite{pu2021distributed,alghunaim2022unified, koloskova2021improved}. In addition, we provide a sufficient condition under which the algorithm achieves linear speedup when the number of iterations is large enough, thereby matching the convergence rates of the centralized stochastic gradient method.

We first present the convergence result of the Stochastic Push-Pull algorithm for minimizing smooth non-convex objective functions, as established in Theorem \ref{thm:convergence}.

\begin{theorem}
    \label{thm:convergence}
    For the Stochastic Push-Pull algorithm outlined in Algorithm \ref{alg:pp}. Let Assumption \ref{a.graph1}, \ref{a.graph}, \ref{a.var} and \ref{a.smooth} hold. Given the number of iterations $T$ with the constants $P_1$ to $P_5$ and $Q$ defined in Lemma \ref{l:descent}, for 
    \[
    \begin{aligned}
        & \gamma = \min \left\{ \prt{\frac{\Delta_f}{\sqrt{\max\crk{P_2,P_3}} \max\crk{N_4, N_7} \sigma^2 L^2}}^{\frac{1}{3}},  \prt{\frac{\Delta_f}{Q\sigma^2 L \prt{T+1}}}^{\frac{1}{2}}, \frac{1}{500\sqrt{\max\crk{P_1,P_2,P_3,P_4}}L} \right\},
    \end{aligned}
    \]
    we have
    \[
    \begin{aligned}
        & \frac{1}{T+1}\sum_{t=0}^{T} \expect\norm{\nabla f(\hat{x}^{(t)})}^2 \lesssim \prt{\frac{Q}{n\pi^2}}^{\frac{1}{2}} \prt{\frac{\Delta_{f}\sigma^2 L}{n(T+1)}}^{\frac{1}{2}}  + \frac{\prt{\Delta_f^2 \sqrt{\max\crk{P_2,P_3}} \max\crk{N_4, N_7} L^2 \sigma^2}^{\frac{1}{3}}}{n\pi}\frac{1}{\prt{T+1}^{\frac{2}{3}}}\\
        &  + \frac{\sqrt{\max\crk{P_1,P_2,P_3, P_4}}L \Delta_f}{n\pi } \frac{1}{T+1}  + \frac{\sqrt{\max\crk{P_2,P_3}} M_1 N_5^2 }{n \pi Q^2} \prt{\frac{\Delta_f L}{\sigma(T+1)}}^2  + \frac{ P_5}{n\pi\prt{T+1}}F_0 .
    \end{aligned}
    \]
\end{theorem}
\begin{proof}
    See Appendix \ref{pf:thm:convergence}.
\end{proof}
For smooth nonconvex objectives, Theorem \ref{thm:convergence} establishes $\mathcal{O}(\frac{1}{\sqrt{T}})$ convergence rate with explicit coefficients for the Push-Pull/AB algorithm under stochastic gradients.

In the following corollary, we discuss the case when $\bR = \bC = \bW$, where $\bW$ is doubly stochastic.  Specifically, we consider three cases: 1) $\bW$ being doubly stochastic only, 2) $\bW$ being symmetric and 3) $\bW$ being positive semi-definite. Case 1 is most general, allowing the communication graph to be either undirected or directed but weight balanced. The other two cases apply to undirected graphs, where Case 3 can be satisfied when choosing large enough diagonal elements in $\bW$.

\begin{corollary}
    \label{co:non-convex}
    For the Stochastic Push-Pull algorithm outlined in Algorithm \ref{alg:pp}. Let Assumption \ref{a.graph1}, \ref{a.graph}, \ref{a.var} and \ref{a.smooth} hold. Denote $\lambda:=\norm{\bW - \mone\mone^\T/n}_2 < 1$. Given the number of iterations $T$, there exists a constant stepsize $\gamma$ such that,\\
    1) \textbf{(doubly stochastic case)} when $\bR = \bC = \bW$, where $\bW$ is doubly stochastic:
    \[
    \begin{aligned}
        & \frac{1}{T+1}\sum_{t=0}^{T}\expect\norm{\nabla f(\hat{x}^{(t)})}^2 \lesssim \prt{\frac{\Delta_f \sigma^2 L}{n(T+1)}}^{\frac{1}{2}}  + \frac{\prt{\Delta_f^2L^2\sigma^2}^{\frac{1}{3}}}{(1-\lambda)\prt{T+1}^{\frac{2}{3}}} \\
        & + \frac{L\Delta_f}{(1-\lambda)^2(T+1)} +  \frac{F_0}{(1-\lambda)^2(T+1)} + \prt{\frac{\sqrt{n}\Delta_f L}{(1-\lambda)^2\sigma(T+1)}}^2 ,
    \end{aligned}
    \]
    2) \textbf{(symmetric case)} when $\bR = \bC = \bW$, where $\bW$ is doubly stochastic and symmetric:
    \[
    \begin{aligned}
        & \frac{1}{T+1}\sum_{t=0}^{T}\expect\norm{\nabla f(\hat{x}^{(t)})}^2 \lesssim \prt{\frac{\Delta_f \sigma^2 L}{n(T+1)}}^{\frac{1}{2}}  + \frac{\prt{\Delta_f^2L^2\sigma^2}^{\frac{1}{3}}}{(1-\lambda)^{\frac{1}{3}}c^{\frac{1}{3}}\prt{T+1}^{\frac{2}{3}}} \\
        & + \frac{L\Delta_f}{\sqrt{c}(1 - \lambda)(T+1)}  + \frac{F_0}{(1-\lambda)(T+1)} + \prt{\frac{\sqrt{n}\Delta_f L }{\sigma(T+1)(1-\lambda)^2}}^2,
    \end{aligned}
    \]
    where $c = 1+\min\crk{\lambda_n,0}$ as in Lemma \ref{l:symmetric}, and\\
    3) \textbf{(positive semi-definite case)} when $\bR = \bC = \bW$, where $\bW$ is doubly stochastic and positive semi-definite:
    \[
    \begin{aligned}
        &\frac{1}{T+1}\sum_{t=0}^{T}\expect\norm{\nabla f(\hat{x}^{(t)})}^2 \lesssim  \prt{\frac{\Delta_f \sigma^2 L}{n(T+1)}}^{\frac{1}{2}}  + \frac{\prt{\Delta_f^2L^2\sigma^2}^{\frac{1}{3}}}{(1-\lambda)^{\frac{1}{3}}\prt{T+1}^{\frac{2}{3}}}  \\
        & + \frac{L\Delta_f}{(1 - \lambda)(T+1)}  + \frac{F_0}{(1-\lambda)(T+1)} + \prt{\frac{\sqrt{n}\Delta_f L }{\sigma(T+1)(1-\lambda)^2}}^2 .
    \end{aligned}
    \]
\end{corollary}
\begin{proof} 
    For a doubly stochastic matrix $\bW$ satisfying Assumption \ref{a.graph}, it holds that $\bpi_\bR = \bpi_\bC = \mone/n$. Consequently, we have $Q = 1/n$ and
    \[
    P_1 = \max\crk{N_8^2, N_5}, P_2 = P_3 =  1, P_4 = N_5, P_5 = N_6/N5. 
    \]
    The results for Case 1, Case 2 and Case 3 directly follow from Lemma \ref{l:symmetric} and Theorem \ref{thm:convergence}. The distinction among these cases arises when estimating $P_5$: in Case 1, note that $N_6 \le N_5^2$, implying $P_5 \le N_5$, whereas in Cases 2 and 3, $P_5$ can be computed exactly. 
\end{proof}
\begin{remark}
\label{remark:case 1}
    In Case 1, we provide a convergence guarantee for the DSGT method under smooth nonconvex objectives when $\bW$ is assumed to be doubly stochastic only. Such a result is new, which particularly demonstrates the linear speedup property of DSGT under such a setting. This finding generalizes upon the previous known results which assume symmetric $\bW$ \cite{alghunaim2022unified,koloskova2021improved}. And the transient time is $\cO(\frac{n^3}{(1-\lambda)^6})$.
\end{remark}
\begin{remark}
\label{remark:case 23}
    In Case 2 and 3, we compare the results in Corollary \ref{co:non-convex} with those obtained for DSGT earlier. Specifically, we calculate and compare the transient times before the term $\cO\prt{\frac{1}{\sqrt{nT}}}$ dominates the remaining terms. In Case 2 (symmetric $\bW$), we obtain a transient time complexity as
    \[
    \cO\prt{\max\crk{\frac{n^3}{c^2(1-\lambda)^2}, \frac{n}{(1-\lambda)^{\frac{8}{3}}}}}.
    \]
    This bound outperforms the best known result in \cite{koloskova2021improved}, which shows $\cO(\max\crk{n/(1-\lambda)^4, n^3/\prt{c^4(1-\lambda)^2}})$. For Case 3 (positive semi-definite $\bW$), we obtain a complexity as
    \[
    \cO\prt{\max\crk{\frac{n^3}{(1-\lambda)^2}, \frac{n}{(1-\lambda)^{\frac{8}{3}}}}},
    \]
    which matches the state-of-the-art result in \cite{alghunaim2022unified}. 
\end{remark}

\subsection{Sufficient Condition for Linear Speedup}
Linear speedup is a desired property of distributed stochastic gradient algorithms, allowing them to match the convergence rate of the centralized SGD algorithm, that is, the $1/\sqrt{nT}$ term can dominate other terms in the complexity bound. However, whether or not Stochastic Push-Pull can achieve linear speedup over general directed graphs has not been discussed in previous works. We address this gap by providing a sufficient condition, stated in Corollary \ref{thm:speedup}, under which the S-PP algorithm achieves linear speedup.
\begin{corollary}
    \label{thm:speedup}
    Let Assumption \ref{a.graph1}, \ref{a.graph}, \ref{a.var} and \ref{a.smooth} hold. Consider the Stochastic Push-Pull algorithm outlined in Algorithm \ref{alg:pp}. Linear speedup is achieved if 
    \[
    \frac{\max\crk{M_1,\tilde{M}_2,M_1 \tilde{M}_2}}{n\prt{\bpi_\bR^\T\bpi_\bC}^2} \le C,
    \]
    where $C$ is a constant that does not depend on $n$.
\end{corollary}
\begin{proof}
    The results follows directly from Theorem \ref{thm:convergence}.
\end{proof}
\begin{remark}
    The sufficient condition applies to many practical settings. For example, it is satisfied when $\bR = \bC = \bW$, where $\bW$ is doubly stochastic as detailed in Corollary \ref{co:non-convex}. 
    
    If only $\bR$ is doubly stochastic (i.e., $\bpi_\bR = \mone/n$), the condition also holds automatically since 
    \[
    \begin{aligned}
        &M_1 = \frac{1}{n}, M_2 = 0, \pi= \frac{1}{n} , \frac{\max\crk{M_1,M_2,M_1 M_2}}{n(\bpi_\bR^\T\bpi_\bC)^2} = 1.
    \end{aligned}
    \] 
    
    Moreover, when $\cG_{\bR}$ and $\cG_{\bC^\T}$ are spanning trees with a common root, matrices $\bR$ and $\bC$ can be chosen to contain only $0$'s and $1$'s and satisfy $M_1 / \prt{n(\bpi_\bR^\T\bpi_\bC)^2} = 1$ \cite{you2025distributed}. Notably, such a design does not use positive diagonals in $\bR$ and $\bC$ as required in previous works such as \cite{pu2020push}.
\end{remark} 

Finally, Corollary \ref{co:linear} indicates that Theorem \ref{thm:convergence} does not allow for higher-order speedup.  
\begin{corollary}
    \label{co:linear}
    Let Assumptions \ref{a.graph1} and \ref{a.graph} hold. We have
    \[
    \frac{\max\crk{M_1,\tilde{M}_2,M_1 \tilde{M}_2}}{n\prt{\bpi_\bR^\T\bpi_\bC}^2} \ge \frac{1}{10}.
    \]
\end{corollary}
\begin{proof}
    When $\cG_{\bR}$ and $\cG_{\bC^\T}$ are spanning trees with a common root, and matrices $\bR$ and $\bC$ are chosen to contain only $0$'s and $1$'s as in \cite{you2025distributed}, we have $M_1 / \prt{n(\bpi_\bR^\T\bpi_\bC)^2} = 1 > \frac{1}{10}$. We then discuss the other cases.
    
    From Assumption \ref{a.graph}, we have $\lim_{k\rightarrow\infty} \bC^{k}=\bpi_\bC\mone^\T$. Then, by the continuity of the norm, we have
    \[
    n\prt{\bpi_\bR^\T\bpi_\bC}^2  = \norm{\lim_{k\rightarrow\infty} \bpi_\bR^\T\bC^{k}}^2 = \lim_{k\rightarrow\infty}\norm{ \bpi_\bR^\T\bC^{k}}^2 .
    \]
    Notice that $\bpi_\bR^\T\bC^{k} = \bpi_\bR^\T\bC + \sum_{t=2}^{k}\bpi_\bR^\T\prt{\bC^t - \bC^{t-1}}$. We have, by picking $\alpha_i = i\norm{\bpi_\bR^\T\prt{\bC^{i+1} - \bC^i}}/M_9$ in Lemma \ref{lem:sum_matrix}, that
    \[
    \begin{aligned}
            & \norm{ \bpi_\bR^\T\bC^{k}}^2 \le 2\norm{ \bpi_\bR^\T\bC }^2 + 2 \norm{\sum_{t=2}^{k}\bpi_\bR^\T\prt{\bC^t - \bC^{t-1}}}^2 \le  2M_1 + 2 M_9 \sum_{t=1}^{k-1} \norm{\bpi_\bR^\T\prt{\bC^{t+1} - \bC^{t}}}/t \\
            \le & 2M_1 + 2 M_9 \sum_{t=1}^{k-1}\prt{\frac{1}{t^2} + \norm{\bpi_\bR^\T\prt{\bC^{t+1} - \bC^{t}}}^2} \le 10 \max\crk{M_1,M_9,M_1 M_9},
    \end{aligned}
    \]
    where we invoked $\sum_{t=1}^{\infty} \frac{1}{t^2} = \pi^2/6\le 2$ in the last inequality.
\end{proof}

\section{Numerical Experiments}
\label{sec:experiment}

In this section, we present numerical experiments to validate the theoretical findings and demonstrate the effectiveness of the Stochastic Push-Pull algorithm in solving decentralized optimization problems. In particular, we compare its performance with other decentralized algorithms, such as SGP \cite{assran2019stochastic} and Push-DIGing \cite{liang2023understanding}, both of which are applicable to arbitrary networks. 

\begin{figure}[htbp]
	\centering
	\subfloat[Erdos-Renyi Graph, $n=8$ and $p=0.1$.]{\includegraphics[width = 0.35\columnwidth]{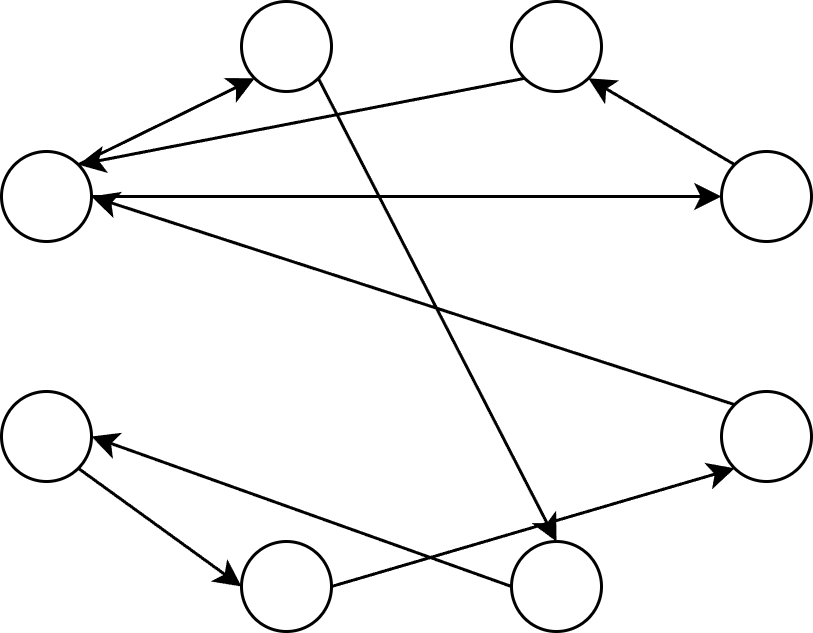}}
    \hfill
	\subfloat[Multi. Sub-Ring Graph, $n = 8$.]{\includegraphics[width = 0.35\columnwidth]{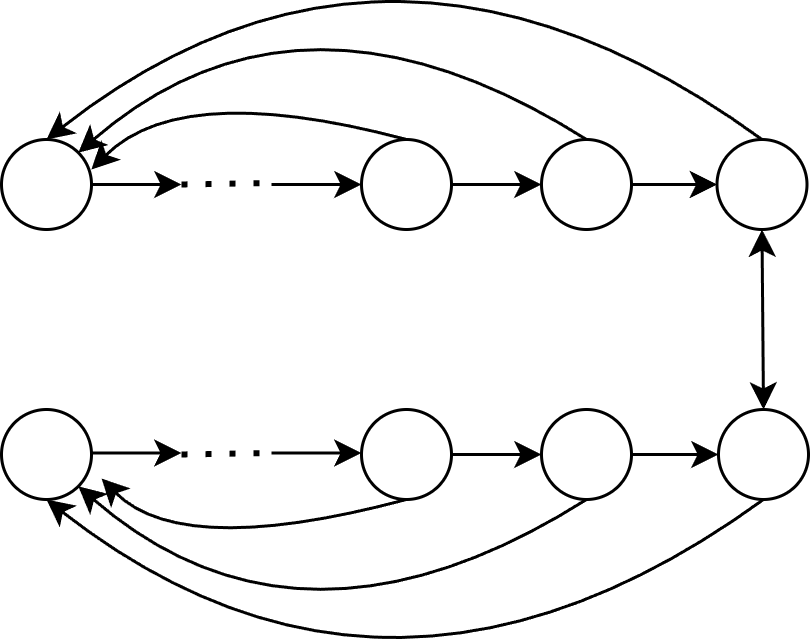}}
	\caption{Illustration of two network topologies. }
	\label{fig:network}
\end{figure}

The network topologies considered in our experiments are illustrated in Figure \ref{fig:network}. The Erdos-Renyi graph (Figure 1(a)) is a directed graph in which each edge is selected with probability $p$. We generate successive Erdos-Renyi graphs until obtaining a strongly connected one. The Multi. Sub-Ring graph (Figure 1(b)) is a directed graph composed of interconnected sub-rings. For simplicity, in S-PP we let $\cG_\bR$ and $\cG_{\bC^\T}$ be the same.

The weight matrices $\bR$ and $\bC$ are constructed by averaging over the in-neighbors and out-neighbors of each node, respectively. Specifically, the entries are given by:
\[
\brk{\bR}_{ij} = \left\{\begin{aligned}
    & \frac{1}{1+d_i^{\text{in}}}, & \text{ if } (j,i)\in \cE_{\bR} \text{ or } i = j,\\
    & 0, & \text{ otherwise},
\end{aligned}\right.
\]
and
\[
\brk{\bC}_{ij} = \left\{\begin{aligned}
    & \frac{1}{1+d_j^{\text{out}}}, & \text{ if } (j,i)\in \cE_{\bC} \text{ or } i = j,\\
    & 0, & \text{ otherwise},
\end{aligned}\right.
\]
where $d_i^{\text{in}}$ and $d_j^{\text{out}}$ are the in-degree of node $i$ and out-degree of node $j$, respectively. Additionally, for SGP and Push-DIGing algorithms,  only the matrix $\bC$ is used as mentioned in \cite{liang2023understanding}. 

\paragraph{Logistic Regression} We compare the performance of S-PP against SGP and Push-DIGing for logistic regression with non-convex regularization \cite{song2022communication,you2024b}. The objective functions $f_i: \mathbb{R}^p\rightarrow \mathbb{R}$ are given by
\[
f_i(x) := \frac{1}{J} \sum_{j = 1}^{J} \ln\prt{1+ \exp(-y_{i,j}h_{i,j}^\T x)} + R\sum_{k = 1}^{p} \frac{x_{[k]}^2}{1+ x_{[k]}^2},
\]
where $x_{[k]}$ is the $k$-th element of $x$, and $\crk{(h_{i,j}, y_{i,j})}_{j=1}^J$ represent the local data kept by node $i$.

To control the data heterogeneity across the nodes, we first let each node $i$ be associated with a local logistic regression model with parameter $\tilde{x}_i$ generated by $\tilde{x}_i = \tilde{x} + v_i$, where $\tilde{x} \sim \cN(0,\bI_p)$ is a common random vector, and $v_i\sim\cN(0, \sigma_h^2\bI_p)$ are random vectors generated independently. Therefore, $\{v_i\}$ decide the dissimilarities between $\tilde{x}_i$, and larger $\sigma_h$ generally amplifies the difference. After fixing $\crk{\tilde{x}_i}$, local data samples are generated that follow distinct distributions. For node $i$, the feature vectors are generated as $h_{i,j} \sim \cN(0, \bI_p)$, and $z_{i,j}\sim\cU(0,1)$. Then, the labels $y_{i,j}\in \crk{-1,1}$ are set to satisfy $z_{i,j}\le 1 + \exp(-y_{i,j}h_{i,j}^\T \tilde{x}_i)$. In the simulations, the parameters are set as follows: $n=20$, $p=400$, $J = 400$, $R = 0.01$, and $\sigma_h = 0.2$. All the algorithms are initialized with the same stepsize $\gamma$ (either 0.1 or 0.05), except for S-PP, which employs a modified stepsize of $\gamma/(n\pi)$ (where $n\pi \approx 1.05$ for Erdos-Renyi Graph and approximately $1.35$ for Multi. Sub-Ring Graph). Such an adjustment is due to S-PP's update mechanism, which incorporates a tracking estimator that effectively accumulates $n\pi$ times the averaged stochastic gradients as the number of iterations increases. Additionally, we implement a stepsize decay of 80\% every $300$ iterations to facilitate convergence, and the experiments are repeated three times to obtain an average performance. 

\begin{figure}
	\centering
	\includegraphics[width=0.7\columnwidth]{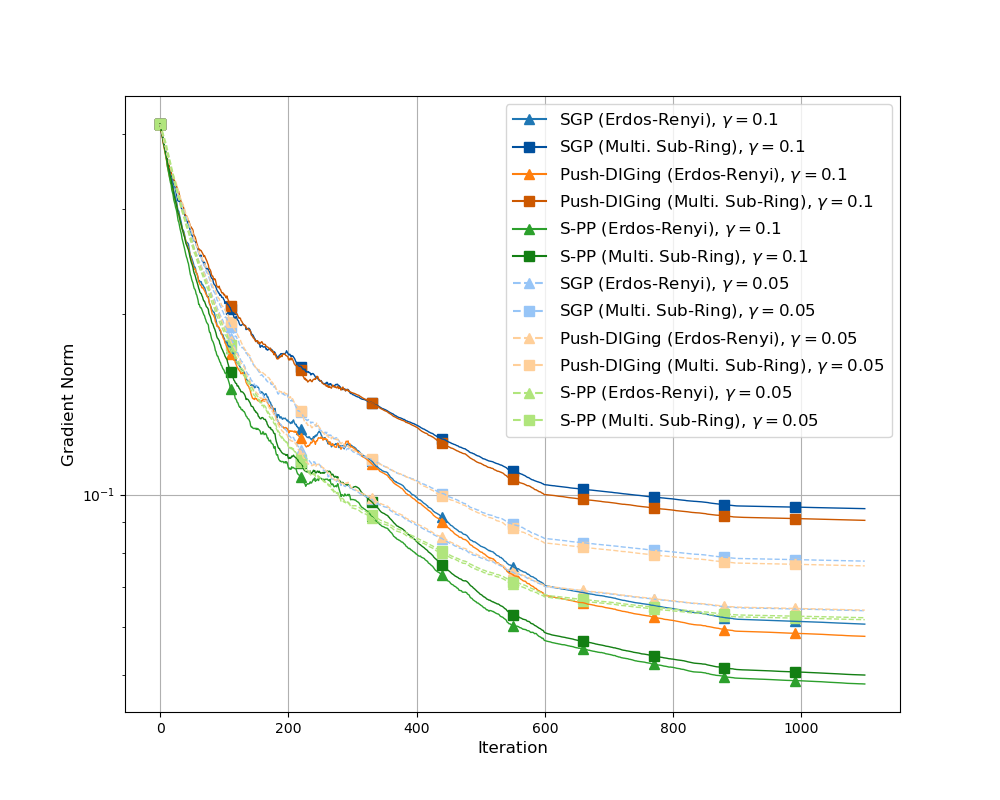}
	\caption{Performance of different algorithms using various communication graphs for logistic regression with nonconvex regularization.}
	\label{fig:logi}
\end{figure}

Figure \ref{fig:logi} illustrates the performance of the algorithms under the two communication graphs. In both network structures, S-PP consistently outperforms the other methods under different stepsize selections. 

\paragraph{Deep Learning} We apply S-PP and the other algorithms to solve the image classification task with CNN over \textbf{MNIST} dataset \cite{lecun2010mnist}. We run all experiments on a server with eight Nvidia RTX 3090 GPUs. The network contains two convolutional layers with max pooling and ReLu and two feed-forward layers. In particular, we consider a heterogeneous data setting, where data samples are sorted based on their labels and then partitioned among the agents. The local batch size is set to $8$ with $20$ agents in total. For all algorithms, the learning rate $\gamma$ is set to 0.01 and decays by 90\% at iterations 8000 and 11000, except for S-PP, which employs a modified stepsize $\gamma/(n\pi)$ for fairness.
\begin{figure}[htbp]
	\centering
	\subfloat{\includegraphics[width = 0.48\columnwidth]{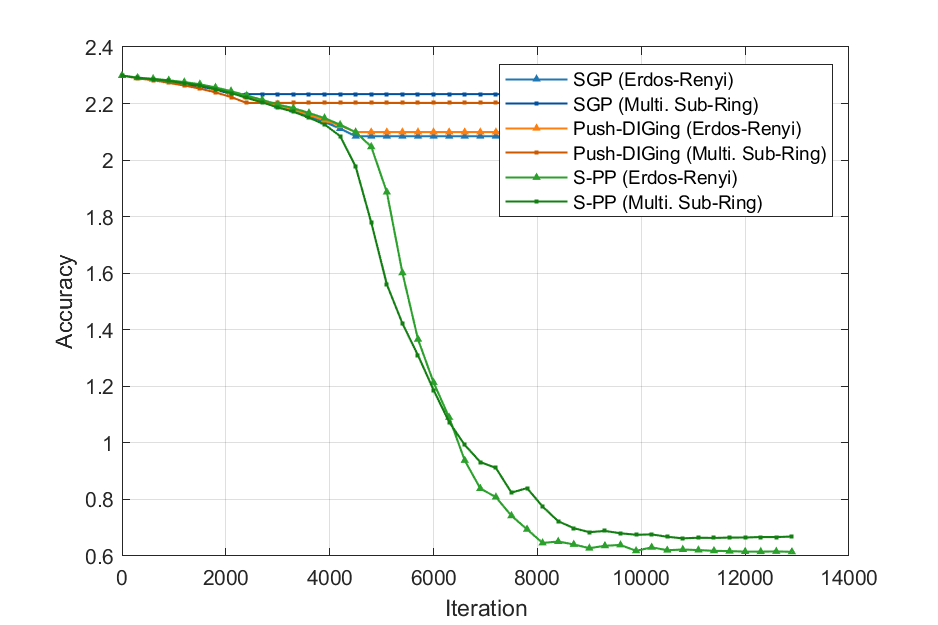}}
    \hfill
	\subfloat{\includegraphics[width = 0.48\columnwidth]{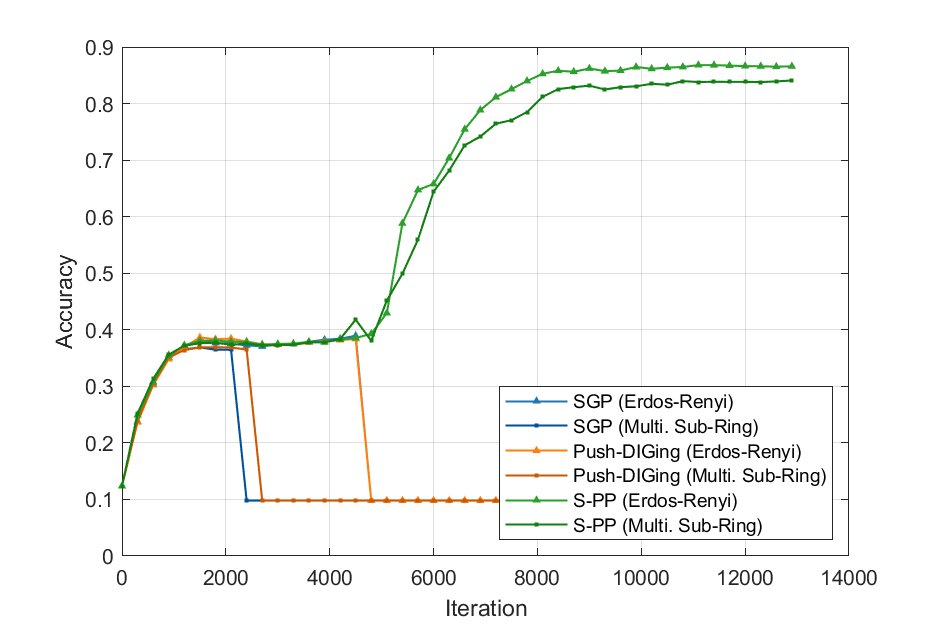}}
	\caption{Training loss and test accuracy of different algorithms for training CNN on MNIST. For the training loss, we plot the curves of the divergent algorithms (SGP (Erdos-Renyi), SGP (Multi. Sub-Ring), Push-DIGing (Erdos-Renyi), and Push-DIGing (Multi. Sub-Ring) with their first non-NaN values in the records.)}
	\label{fig:mnist}
\end{figure}

Figure \ref{fig:mnist} displays the training loss and test accuracy curves, indicating that S-PP outperforms the other methods. In contrast, SGP and Push-DIGing are notably sensitive to both the stepsize and the network size $n$. In particular, they may fail to converge when $n$ is large or the step-size is chosen too aggressively.

\section{Conclusions}
\label{sec:conclusion}
In this paper, we present a comprehensive convergence analysis for the Stochastic Push–Pull algorithm, with particular emphasis on the underlying assumptions and convergence rate. The proposed analytical framework offers a general approach to studying decentralized optimization algorithms and encompasses a wide range of state-of-the-art results on decentralized stochastic gradient tracking methods. Furthermore, we establish a sufficient condition under which the Stochastic Push–Pull algorithm achieves linear speedup. Our theoretical findings are validated through numerical experiments, which demonstrate the effectiveness of the Stochastic Push–Pull algorithm in solving decentralized optimization problems over directed networks.

\appendix

\section{Proofs}

\subsection{Proof of Lemma \ref{l:eigen}}

\label{pf:l:eigen}

We prove the lemma by first examining the structure of the graph $\cG_{\bR}$. 

If $\cR_{\bR} = [n]$, meaning that there is at least one spanning tree enabling the existence of a path between any two nodes, then $\cG_{\bR}$ is strongly connected. By the Perron-Frobenius theorem, there exists a unique positive eigenvector $\bpi_{\bR}$, corresponding to the eigenvalue $1$ satisfying $ \bpi_{\bR}^\T\bR = \bpi_{\bR}^\T$ and $\norm{\bpi_\bR}_1 = 1$. Consequently, this guarantees the existence and uniqueness of the root eigenvector for $\bR$.

If $\cR_{\bR} \ne [n]$, there exists $n_0 < n$ such that $|\cR_{\bR}| = n_0$. Thus, there is no direct path from any node $i\notin \cR_{\bR}$ to any root nodes $j\in \cR_{\bR}$. In this case, we can reorder the vertices such that $\bR$ is represented as $\hat{\bR} = \prt{\begin{matrix}
        \bR_1 & \mathbf{0}\\
       \bR_2 & \bR_3  
\end{matrix}} $ where $\bR_1\in \reals^{n_0\times n_0}$ corresponds exclusively to the root nodes in $\cR_{\bR}$. Here, $\cG_{\bR_1}$ is strongly connected, and by the Perron-Frobenius theorem, $\bR_1$ has a unique positive left eigenvector $\bpi_{\bR_1}$ corresponding to the eigenvalue $1$ (where $1$ is the spectral radius). Suppose there exist two root eigenvectors $\brk{\mathbf{u}_1^\T, \mathbf{0}^\T}^\T$ and $\brk{\mathbf{u}_2^\T, \mathbf{0}^\T }^\T$ associated with $\hat{\bR}$; then, we must have $\bu_1 = \bu_2 = \bpi_{\bR_1}$. Thus, the unique root eigenvector of $\bR$ can be reordered as $\brk{\bpi_{\bR_1}^\T, \mathbf{0}^\T}^\T$.

Similarly, a unique root eigenvector exists for $\bC^\T$. Furthermore, since $\cR := \cR_{\bR}\cap \cR_{\bC} \ne \emptyset$, we have
  \[
  \bpi_\bR^\T\bpi_\bC = \sum_{i\in \cR} \bpi_{\bR}(i) \bpi_{\bC}(i) > 0.
  \]

\subsection{Proof of Lemma \ref{lem:diagonal}}
\label{pf:lem:diagonal}

We begin by introducing the spectral radius.
\begin{definition}{(spectral radius)}
    \label{def:spectral_radius}
    Let $\lambda_1,\cdots, \lambda_n$ be the eigenvalues of a matrix $\bA\in \mathbb{C}^{n\times n}$. Then, its spectral radius $\rho(\bA)$ is defined as:
    \[
    \rho(\bA) := \max\crk{|\lambda_1|, \cdots, |\lambda_n|}.
    \]
\end{definition}

\begin{lemma}
    \label{l:spectral_radius}
    Given any row-stochastic and nonnegative matrix $\bW\in \reals^{n\times n}$ with strictly positive diagonal. If the induced graph $\cG_\bW$ contains at least one spanning tree, the following holds.
    \begin{enumerate}
        \item The spectral radius $\rho(\bW) = 1$ and $1$ is an eigenvalue of $\bW$ corresponding to the eigenvector $\mone$.
        \item There exists the root eigenvector $\bpi\in \reals^{n}$ of eigenvalue $1$ such that $\bpi^\T\bW = \bpi^\T$ and $\bpi^\T\mone = 1$.
        \item The eigenvalue $1$ of $\bW$ has algebraic multiplicity equal to 1, and $1$ is the unique eigenvalue of maximum modulus.
    \end{enumerate}
\end{lemma}
\begin{proof}
    Statements 1 and 3 follow from \cite[Lemma 3.4]{ren2005consensus}. Statement 2 follows from choosing $\bR = \bW$ in Lemma \ref{l:eigen}.
\end{proof}
\begin{lemma}{(Gelfand's formula, see \cite[Corollary 5.6.14]{horn2012matrix})}
    \label{l:gelfand_formula}
    For any matrix norm $\norm{\cdot}$, we have
    $
    \rho(\bA) = \lim_{k\rightarrow \infty} \norm{\bA^k}^{1/k}.
    $
\end{lemma}

We then prove Lemma \ref{lem:diagonal} as follows.
    First, let $\crk{\lambda_i}$ and $\crk{\bv_i}$ be the eigenvalues and corresponding right eigenvectors of $\bA$, where $\lambda_1 = 1$ and $\bv_1 = \mone$. By Lemma \ref{l:spectral_radius}, we have $|\lambda_i| < 1$ for $i\neq 1$. Then,
    \[
    \begin{aligned}
        \prt{\bA - \mone\bpi^\T}\bv_1 & = \prt{\bA - \mone\bpi^\T}\mone = 0, \\
        \prt{\bA - \mone\bpi^\T}\bv_i & = \lambda_i\bv_i - \mone\bpi^\T\bv_i = \lambda_i\bv_i , \qquad \forall \ i >1,
    \end{aligned}
    \]
    where we use $\bpi^\T\bv_i = 0$ in the second line, since
    \[
    (1-\lambda_i)\bpi^\T\bv_i = \bpi^\T\bv_i - \lambda_i\bpi^\T\bv_i = \bpi^\T\bA\bv_i - \bpi^\T \bA\bv_i = 0.
    \]
    Hence, the spectrum of $\bA - \mone\bpi^\T$ are 
    $
    \crk{0, \lambda_2, \cdots, \lambda_n},
    $
    and thus the spectral radius $\rho(\bA - \mone\bpi^\T) = \lambda := \max\crk{|\lambda_2|, \cdots, |\lambda_n|} < 1$. 
    
    Next, observe that 
    \[
    \bA^{m} - \mone\bpi^\T = \prt{\bA - \mone\bpi^\T}^m.
    \]
    Applying Gelfand's formula (Lemma \ref{l:gelfand_formula}) to $\bA - \mone\bpi^\T$ with $2$-norm, we have for any $\epsilon \in \prt{0, 1- \rho(\bA - \mone\bpi^\T)}$, there exists an integer $m_{\bA,\epsilon}>0$ such that for any $m\ge m_{\bA,\epsilon}$,
    \[
    \begin{aligned}
    & \norm{\bA^m - \mone\bpi^\T}_2 = \norm{\prt{\bA - \mone\bpi^\T}^m}_2 \\
    \le & \prt{\rho(\bA - \mone\bpi^\T) + \epsilon}^m = \prt{\lambda + \epsilon}^m < 1.
    \end{aligned}
    \]
    Choosing $\epsilon = p = \frac{1}{2}\prt{1- \lambda }$ and setting $m_{\bA} =  m_{\bA,\epsilon}$ completes the proof.

\subsection{Proof of Lemma \ref{l:symmetric}}
\label{pf:l:symmetric}
Define $\bPi := \bI - \frac{1}{n}\mone\mone^\T$, $\tbW = \bPi\bW$ and $\lambda := \norm{\bPi\bW}_2 = \max\crk{|\lambda_2|, |\lambda_n|}\in (0,1)$. It follows that $\bpi_\bR = \bpi_\bC = \mone/n$. Then,
\[
M_1 = 1/n, \quad M_2 = N_1 = N_2 = N_3 = N_4 = 0.
\]
With $\lambda = \norm{\bW - \mone\mone^\T/n}_2 < 1$, we have
\[
\begin{aligned}
    N_5 = & \sum_{t=1}^{\infty} \norm{(t-1)\tbW^{t}}_2\le \sum_{t=1}^{\infty} (t-1)\lambda^t \le \frac{\lambda^2}{(1-\lambda)^2}, N_6  =  \sum_{t=1}^{\infty}\norm{t\tbW^t}_2^2 \le \sum_{t=1}^{\infty} t^2\lambda^{2t}\le \frac{\lambda^2 + \lambda^4}{(1-\lambda^2)^3},\\
    N_7 \le & 2\sum_{t=1}^{\infty} \norm{t\tilde{\bW}^{t+1}}_2^2 + 2\sum_{t=1}^{\infty} \norm{(t-1)\tilde{\bW}^{t}}_2^2 \le \frac{8}{(1-\lambda^2)^3}, N_8 \le  \sum_{t=1}^{\infty} \norm{t\tilde{\bW}^{t+1}}_2 + \sum_{t=1}^{\infty} \norm{(t-1)\tilde{\bW}^{t}}_2 \le \frac{8}{(1-\lambda)^2}.
\end{aligned}
\]
If $\bW$ is additionally symmetric, let $\lambda_1,\cdots,\lambda_n$ be the eigenvalues of $\bW$ with $\lambda_1 = 1$. By the Perron-Frobenius theorem and Assumption \ref{a.graph}, $|\lambda_i| < 1$ for all $i\ne 1$. Without loss of generality, assume that $\lambda_1 = 1>\lambda_2\ge\lambda_3\ge \cdots\ge \lambda_n > -1$. It follows directly from the definition of the symmetric matrix $\bR = \bC = \bW$ that
\[
\begin{aligned}
    M_1 & = 1/n, \quad M_2 = N_1 = N_2 = N_3 = N_4 = 0,N_5  = \sum_{t=1}^{\infty}\norm{(t-1)\tbW^t}_2 , N_6 = \sum_{t=1}^{\infty}\norm{t\tbW^t}_2^2,\\
    N_7 & = \sum_{t=1}^{\infty}\norm{t\tbW^{t+1} - (t-1)\tbW^{t}}_2^2, N_8 = \sum_{t=1}^{\infty}\norm{t\tbW^{t+1} - (t-1)\tbW^{t}}_2.
\end{aligned}
\]
Note that the eigenvalues of $\tbW^t$ for any $t$ are $\crk{0,\lambda_2^t,\cdots,\lambda_n^t}$ and $\norm{\tbW^t}_2 = \max\crk{|\lambda_2|^t, |\lambda_n|^t} = \lambda^t$, implying that
\[
N_5 = \sum_{t=1}^{\infty}(t-1)\lambda^t = \frac{\lambda^2}{(1-\lambda)^2},  N_6 = \sum_{t=1}^{\infty}t^2\lambda^{2t} = \frac{\lambda^2+\lambda^4}{(1-\lambda^2)^3}.
\]
We illustrate the estimation of $N_7$ and $N_8$ with a useful inequality in Lemma \ref{lem:W-series}.
\begin{lemma}
    \label{lem:W-series}
    Let $\lambda \in (-1,1)$. Then, $|t\lambda^t|\le \frac{1}{1-|\lambda|}$ for all $t\ge 0$.
\end{lemma}
\begin{proof}
We have
    \[
    |t\lambda^t|\le t|\lambda|^t \le \sum_{j=1}^{t}|\lambda|^j \le \frac{1 - |\lambda|^{t+1}}{1-|\lambda|} \le \frac{1}{1-|\lambda|}.
    \]
\end{proof}

Note that all the eigenvalues of $t\tbW^{t+1} - (t-1)\tbW^{t}$ are $f_t(x)$ defined as $tx^{t+1} - (t-1)x^{t}$ where $x\in \Lambda := \crk{0,\lambda_2,\cdots,\lambda_n}$. Then,
\[
\norm{t\tbW^{t+1} - (t-1)\tbW^{t}}_2^2 = \max_{x\in \Lambda} \crk{|f_t(x)|}^2 =  \max_{x\in \Lambda} \crk{\prt{f_t(x)}^2}.
\]
If the maximum is obtained by a positive $\hat{x}\in \Lambda$, we have
\[
\begin{aligned}
    &\prt{f_t(\hat{x})}^2  = \prt{t\hat{x}^{t+1} - (t-1)\hat{x}^t}^2 \\
    =& \prt{\hat{x}^{t+1} - (t-1)\hat{x}^t(1-\hat{x})}^2 \\
     \le & 2\prt{\hat{x}}^{2t+2} + 2(t-1)^2 \hat{x}^{2t}(1-\hat{x})^2\\
    \le & 2\prt{\hat{x}}^{2t+2} +  2\hat{x}^{t}\prt{1+\sqrt{\hat{x}}}^2 \le  10 \hat{x}^t,
\end{aligned}
\]
where the second inequality comes from Lemma \ref{lem:W-series}, i.e., $(t-1)\sqrt{\hat{x}}^t \le 1/(1-\sqrt{\hat{x}} )$, and the last inequality holds by the fact that $\hat{x}\in (0,1)$.

If the maximum is obtained by a negative $\hat{x}\in \Lambda$, we have
\[
\begin{aligned}
    & \prt{f_t(\hat{x})}^2  = \prt{t\hat{x}^{t+1} - (t-1)\hat{x}^t}^2 \\
    \le & 2t^2 \hat{x}^{2t+2} + 2(t-1)^2\hat{x}^{2t} \le 2\prt{t^2|\hat{x}|^{t}+ (t-1)^2|\hat{x}|^{t-1}}|\hat{x}|^{t} \\
    \le & \frac{4}{\prt{1-\sqrt{|\hat{x}|}}^2}|\hat{x}|^{t} \le \frac{4}{1-|\hat{x}|}|\hat{x}|^{t} \le \frac{4}{1-|\lambda_n|}|\hat{x}|^{t},
\end{aligned}
\]
where the third inequality is from Lemma \ref{lem:W-series}, and the last inequality holds because $\lambda_n \le \hat{x} <0$ given $\hat{x}\in \Lambda$ and $\hat{x}<0$.

Since $\hat{x}\in \Lambda$, it holds that $|\hat{x}| \le \lambda$. Combining the two cases, we have, for any $t\ge 1$,
\[
\begin{aligned}
    & \norm{t\tbW^{t+1} - (t-1)\tbW^{t}}_2^2 \le \frac{10}{1 + \min\crk{\lambda_n,0}}\lambda^t, \\
    & \norm{t\tbW^{t+1} - (t-1)\tbW^{t}}_2  \le \frac{\sqrt{10}}{\sqrt{1 + \min\crk{\lambda_n,0}}}\lambda^{\frac{t}{2}}.
\end{aligned}
\]
Therefore, we conclude that
\[
\begin{aligned}
    N_7  & \le  \sum_{t=1}^{\infty} \frac{10}{1 + \min\crk{\lambda_n,0}}\lambda^t \le \frac{10}{\prt{1 + \min\crk{\lambda_n,0}}\prt{1-\lambda}} \\
    N_8 & \le \sum_{t=1}^{\infty} \frac{\sqrt{10}}{\sqrt{1 + \min\crk{\lambda_n,0}}}\lambda^{\frac{t}{2}} \le \frac{10}{\sqrt{1 + \min\crk{\lambda_n,0}}(1 - \lambda)} ,\\
\end{aligned}
\]
where we use the fact $\frac{1}{1 - \sqrt{\lambda}}  \le \frac{2}{1 - \lambda}$ in the last inequality.

\subsection{Proof of Lemma \ref{lem:series_help1}}
\label{pf:lem:series_help1}
    Note that $\sum_{t=k+1}^{T}\sum_{j=1}^k a_{t-j}$ only consists of the terms $a_i$ for $i = 1,\cdots, T-1$. We denote
    \[
    \sum_{t=k+1}^{T}\sum_{j=1}^k a_{t-j} := \sum_{i=1}^{T-1} a_i b_i,
    \]
    and show that $b_i \le i$.
    
    Define the set $S_{t,j,k}:= \crk{\prt{t,j}: k+1\le t\le T, 1\le j\le k}$. By the closed form of $b_i$, we have for any given $i$ ($1\le i \le T-1$) that
    \[
    \begin{aligned}
        & b_i = \sum_{t=k+1}^{T}\sum_{j=1}^k \mone_{\crk{t-j=i} } =  \sum_{\prt{t,j}\in S_{t,j,k}}  \mone_{\crk{t-j=i} } = \sum_{\prt{t,j}\in S_{t,j,k}\cap \crk{t-j=i}} 1.
    \end{aligned}
    \]
    Denote $\tilde{S}_{t,j,k,i}:= S_{t,j,k}\cap \crk{t-j=i}$ given $i,k$. Then,
    \[
    \begin{aligned}
        & \tilde{S}_{t,j,k,i}  = \crk{\prt{t,j}: k+1\le t\le T, 1\le j\le k, t-j = i} \\
        = & \crk{\prt{t,j}: \max\crk{1,k+1-i}\le j\le \min\crk{k,T-i}, t = j+i} \\
        \subseteq & \crk{\prt{t,j}: k+1-i\le j\le k, t = j+i}.
    \end{aligned}
    \]
    Thus, we have 
    \[
    b_i = |\tilde{S}_{t,j,k,i}| \le |\crk{\prt{t,j}: k+1-i\le j\le k, t = j+i}| = i ,
    \]
    which completes the proof.

\subsection{Proof of Lemma \ref{l:variance}}
\label{pf:l:variance}
Note that for any given $t\ge 0$, each sample $\xi_i^{(t)}$ is drawn independently from the distribution $\cD_{i}$ for all $i\in [n]$ by Assumption \ref{a.var}. Then, we have, for any $i\ne j \in [n]$, that
\[
\begin{aligned}
    & \expect\left\langle \theta_i^{(t)}, \theta_j^{(t)} \right\rangle =  \expect \brk{ \expect \left\langle \theta_i^{(t)} , \theta_j^{(t)} \right\rangle | \cF_t, \xi_i^{(t)}} = \expect \brk{ \left\langle  \theta_i^{(t)}, \expect \theta_j^{(t)}| \cF_t, \xi_i^{(t)} \right\rangle  } =  \expect \brk{ \left\langle  \theta_i^{(t)}, \expect \theta_j^{(t)} | \cF_t  \right\rangle  } =  0.
\end{aligned}
\]
Thus, under Assumption \ref{a.var}, we have that
\[
\expect\brk{\norm{a^\T \btheta^{(t)}}^2} = \sum_{i=1}^{n} \expect\brk{a_i^2 \norm{\theta_i^{(t)}}^2} \le  \norm{a}^2\sigma^2.
\]

\subsection{Proof of Lemma \ref{l:X_diff}}
\label{pf:l:X_diff}
Notice that
    \[
    \begin{aligned}
       &  \bar{\bX}^{(t+1)} - \bar{\bX}^{(t)} = - \gamma \mone \bpi_\bR^\T \bY^{(t)}  =  - \gamma \mone \bpi_\bR^\T  \prt{\bPi_{\bC}\bY^{(t)} + \bpi_{\bC}\mone^\T \bY^{(t)}}.\\
    \end{aligned}
    \]
    Multiplying $\mone^\T$ on both sides of Equation \eqref{eq:pp2}, we have $\mone^\T \bY^{(t)} = \mone^\T\bG^{(t)}$.
    
    Then, implementing the Equation \eqref{eq:pp2_PiC} yields
    \[
        \begin{aligned}
            \Delta \bar{\bX}^{(t)}  = & - \gamma \mone \bpi_\bR^\T \sum_{j=1}^{t-1} \prt{\tbC^{t-j+1} - \tbC^{t-j}} \bG^{(j)}  - \gamma \mone \bpi_\bR^\T  \tbC \bG^{(t)} - \gamma \mone \bpi \mone^\T \bG^{(t)}.
        \end{aligned}
    \]
    Consequently, it holds that
    \[
    \begin{aligned}
        & \Delta\bar{\bX}^{(t)}  =  - \gamma \mone  \sum_{j=1}^{t-1}\bpi_\bR^\T \prt{\bC^{t-j+1} - \bC^{t-j}}\btheta^{(j)} - \gamma \mone \bpi_\bR^\T \bC \btheta^{(t)}  \\
        & - \gamma \mone \bpi_\bR^\T \sum_{j=1}^{t-1} \prt{\tbC^{t-j+1} - \tbC^{t-j}} \nabla \bF^{(j)}  - \gamma \mone \bpi_\bR^\T \bC \nabla \bF^{(t)} \\
        & :=  \bQ_{1,t} + \bQ_{2,t} ,
    \end{aligned}
    \]
    where we denote $\bQ_{1,t}, \bQ_{2,t}$ as the terms related to $\btheta$ and $\nabla \bF$ in the decomposition of $\Delta\bar{\bX}^{(t)}$ respectively. Next, we show that the components in $\bQ_{1,t}$ and $\bQ_{2,t}$ are upper bounded in the following two steps, respectively.
    
    \textbf{Step 1}: Taking full expectation on the squared F-norm of $\bQ_{1,t}$, invoking Assumption \ref{a.var} and Lemma \ref{l:co-var}, we have
    \[
    \begin{aligned}
         & \expect \norm{\bQ_{1,t}}_F^2 \\
        = & \gamma^2 n\expect \norm{\sum_{j=1}^{t-1}\bpi_\bR^\T \prt{\bC^{t-j+1} - \bC^{t-j}} \btheta^{(j)} + \bpi_\bR^\T \bC \btheta^{(t)}}_F^2 \\
        \le & \gamma^2 n \sum_{j=1}^{t-1}  \norm{ \bpi_\bR^\T \prt{\bC^{t-j+1} - \bC^{t-j}} }_2^2 \sigma^2 + \gamma^2 n \norm{\bpi_\bR^\T \bC}_2^2 \sigma^2 \\
        \le & \gamma^2 n M_1 \sigma^2.
    \end{aligned} 
    \]
    \textbf{Step 2}:     Inspired by the transformation in Equation \eqref{eq:pp2_PiC}, we have
    \[
    \begin{aligned}
        \bQ_{2,t}  = & - \gamma \mone \bpi_\bR^\T \sum_{j=1}^{t-1}\prt{\tbC^{t-j+1} - \tbC^{t-j}} \nabla \bF^{(j)}  - \gamma \mone \bpi_\bR^\T \tbC \nabla \bF^{(t)} - \gamma \mone \bpi_\bR^\T \bpi_\bC \mone^\T \nabla \bF^{(t)} \\
        = & -\gamma \mone\bpi_\bR^\T \sum_{j=0}^{t-1} \tbC^{t-j} \brk{ \nabla \bF^{(j+1)} - \nabla \bF^{(j)}}  - \gamma \mone\bpi_\bR^\T \tbC^{t} \nabla \bF^{(0)} - \gamma \mone \bpi_\bR^\T \bpi_\bC \mone^\T \nabla \bF^{(t)}.
    \end{aligned}
    \]
    Taking the square of F-norm on $\bQ_{2,t}$ and then the full expectation, we have
    \[
    \begin{aligned}
        & \expect\norm{\bQ_{2,t   } }_F^2 \le  3 \gamma^2 n \expect \norm{\sum_{j=0}^{t-1}\bpi_\bR^\T  \tbC^{t-j} \brk{ \nabla \bF^{(j+1)} - \nabla \bF^{(j)} } }_F^2 \\
        &  + 3 \gamma^2 n \expect\norm{\pi \mone^\T \nabla \bF^{(t)}}_F^2 + 3 \gamma^2 n \expect\norm{ \bpi_\bR^\T  \tbC^{t} \nabla \bF^{(0)} }_F^2.
    \end{aligned}
    \]
    To further bound $\expect\norm{\bQ_{2,t   } }_F^2$, firstly, invoking Assumption \ref{a.smooth}, Lemma \ref{lem:matrix_norm} and Lemma \ref{lem:sum_matrix} (where we choose $\alpha_{j} = \norm{ \bpi_\bR^\T \tbC^{t-j}} /N_1$ for $0\le j\le t-1$), 
    \[
    \begin{aligned}
        & \expect \norm{\sum_{j=0}^{t-1}\bpi_\bR^\T  \tbC^{t-j} \brk{ \nabla \bF^{(j+1)} - \nabla \bF^{(j)}}}_F^2 \\
        \le & \sum_{j=0}^{t-1} \frac{1}{\alpha_j} \expect \norm{\bpi_\bR^\T  \tbC^{t-j} }_2^2 \norm{ \nabla \bF^{(j+1)} - \nabla \bF^{(j)}  }_F^2 \\
        \le & 3N_1 L^2 \sum_{j=0}^{t-1} \norm{\bpi_\bR^\T \tbC^{t-j} }  \expect \brk{\norm{\hat{ \bX}^{(j+1)}}_F^2 + \norm{\hat{\bX}^{(j)}}_F^2 + \norm{\Delta\bar{\bX}^{(j)}}_F^2} .
    \end{aligned}
    \]
    Secondly, it follows from Assumption \ref{a.smooth} that
    \[
    \begin{aligned}
        & \expect\norm{\pi \mone^\T \nabla \bF^{(t)} }^2  \le  2\pi^2 \expect \brk{\norm{ \mone^\T  \nabla \bbF^{(t)} }^2 +  \norm{  \mone^\T \nabla \bF(\bar{\bX}^{(t)}) }^2} \\
        & \le  2 n\pi^2 L^2\expect\norm{\hat{\bX}^{(t)}}_F^2  + 2  n^2\pi^2 \expect\norm{\nabla f(\hat{x}^{(t)})}^2
    \end{aligned}
    \]
    Finally, by Lemma \ref{lem:matrix_norm},
    \[
    \expect\norm{ \bpi_\bR^\T \tbC^{t} \nabla \bF^{(0)} }_F^2\le \norm{\bpi_\bR^\T \tbC^{t}  }^2\expect \norm{\nabla \bF^{(0)} }_F^2.
    \]
    
    Note that $\expect \norm{\Delta\bar{\bX}^{(t)}  }_F^2  \le 2\expect \norm{\bQ_{1,t}}_F^2 + 2\expect \norm{\bQ_{2,t}}_F^2$.
     We take the summation over $t$ from $0$ to $T$ for each term in the bounds of $\expect \norm{\bQ_{1,t}}_F^2$ and $\expect \norm{\bQ_{2,t}}_F^2$.
    Invoking Lemma \ref{lem:sum_help}, letting $a_i = \norm{\bpi_\bR^\T  \bPi_\bC \bC^{i} } $ and $b_i = \norm{\hat{\bX}^{(i+1)}}_F^2$, $\norm{\hat{\bX}^{(i)}}_F^2$, $ \norm{\Delta\bar{\bX}^{(i)} }_F^2$ respectively, we get, for instance, 
    \[
    \sum_{t=0}^{T} \sum_{j=0}^{t-1} \expect\norm{\bpi_\bR^\T  \tbC^{t-j} } \norm{\hat{\bX}^{(j+1)}}_F^2 \le N_1\sum_{t=0}^{T}\expect\norm{\hat{\bX}^{(t)}}_F^2.
    \]
     Combining all such inequalities yields
    \[
    \begin{aligned}
        & \sum_{t=0}^{T} \expect \norm{\Delta\bar{\bX}^{(t)} }_F^2  \le 2\gamma^2 n M_1 \sigma^2 (T+1) \\
        & + 48 \gamma^2 n \max\crk{N_1^2, n \pi^2} L^2 \sum_{t=0}^{T} \expect\norm{\hat{\bX}^{(t)}}_F^2 \\
        &  + 18 \gamma^2 n N_1^2 L^2 \sum_{t=0}^{T} \expect\norm{\Delta\bar{\bX}^{(t)}  }_F^2 \\
        & + 12 \gamma^2  n^3 \pi^2\sum_{t=0}^{T} \expect\norm{\nabla f(\hat{x}^{(t)})}_2^2   + 6 \gamma^2 n N_2  \norm{\nabla \bF^{(0)}}_F^2.
    \end{aligned}
    \]
    Thus, for $\gamma \le \frac{1}{6\sqrt{n} N_1 L}$ (hence $18 \gamma^2 n N_1^2 L^2 \le \frac{1}{2}$), we obtain the desired result by rearranging the terms.

\subsection{Proof of Lemma \ref{l:PiX}}
\label{pf:l:PiX}

Multiplying both sides of Equation \eqref{eq:pp1} by $\bPi_{\bR}$, we obtain
\begin{equation}
    \label{eq:PiX-0}
    \begin{aligned}
        \hat{\bX }^{(t)} = & - \gamma \sum_{j=0}^{t-2}\sum_{k=1}^{t-j-1}\tbR^{k} \tbC^{t-j-k} \brk{\bG^{(j+1)} - \bG^{(j)}} \\
        &- \gamma \sum_{j=0}^{t-2}\sum_{k=1}^{t-j-1}\tbR^{k} \bpi_{\bC}\mone^\T \brk{\bG^{(j+1)} - \bG^{(j)}} \\
        &  - \gamma \sum_{k=1}^{t}  \tbR^{k} \tbC^{t-k} \bG^{(0)} - \gamma \sum_{k=1}^{t}  \tbR^{k} \bpi_\bC\mone^\T \bG^{(0)} \\
        := &  - \gamma \bH_{1,t} - \gamma \bH_{2,t}.
    \end{aligned}
\end{equation}
We bound $\bH_{1,t}$ and $\bH_{2,t}$ by rearranging the terms in Equation \eqref{eq:PiX-0},  as outlined in the following steps.

\textbf{Step 1}: For $\bH_{1,t}:= \sum_{j=0}^{t-2}\sum_{k=1}^{t-j-1}\tbR^{k} \bpi_{\bC}\mone^\T \brk{\bG^{(j+1)} - \bG^{(j)}}+ \sum_{k=1}^{t}  \tbR^{k} \bpi_\bC\mone^\T \bG^{(0)}$, it holds that
\[
\begin{aligned}
    &\bH_{1,t} =  \sum_{j=0}^{t-1} \tbR^{t-j} \bpi_\bC\mone^\T \btheta^{(j)} + \sum_{j=0}^{t-1} \tbR^{t-j} \bpi_\bC\mone^\T \nabla \bF^{(j)}.
\end{aligned}
\]
Invoking Assumption \ref{a.var} and Lemma \ref{l:co-var}, we have
\[
\begin{aligned}
    &\expect \norm{\sum_{j=0}^{t-1}\tbR^{t-j} \bpi_\bC\mone^\T \btheta^{(j)} }_F^2 = \sum_{j=0}^{t-1} \expect \norm{ \tbR^{t-j} \bpi_\bC \mone^\T \btheta^{(j)}}_F^2 \\
    \le & \sum_{j=0}^{t-1}\norm{ \tbR^{t-j} \bpi_\bC}_2^2 n\sigma^2  \le n N_4 \sigma^2. \\
\end{aligned}
\]
Furthermore, observe that
\[
\begin{aligned}
    & \sum_{j=0}^{t-1} \tbR^{t-j} \bpi_\bC\mone^\T \nabla \bF^{(j)}
    = \sum_{j=0}^{t-1} \tbR^{t-j} \bpi_\bC\mone^\T \nabla\bbF^{(j)}  \\
    & \qquad + \sum_{j=0}^{t-1} \tbR^{t-j} \bpi_\bC n \nabla f(\hat{x}^{(j)}) .\\
\end{aligned}
\]
By choosing $\alpha_{t-j} = \norm{ \tbR^{t-j} \bpi_\bC}/N_3$ in Lemma \ref{lem:sum_matrix},  we obtain the following bounds:
\[
\begin{aligned}
    &  \norm{\sum_{j=0}^{t-1} \tbR^{t-j} \bpi_\bC\mone^\T  \nabla\bbF^{(j)} }_F^2 \le N_3 \sum_{j=0}^{t-1}  \norm{ \tbR^{t-j} \bpi_\bC} n \norm{\nabla\bbF^{(j)}}_F^2 \\
    & \le  N_3 n L^2\sum_{j=0}^{t-1}  \norm{ \tbR^{t-j} \bpi_\bC} \norm{\hat{\bX}^{(j)}}_F^2,\\
    &\norm{ \sum_{j=0}^{t-1}\tbR^{t-j} \bpi_\bC n \nabla f(\hat{x}^{(j)}) }_F^2  \le  N_3 n^2 \sum_{j=0}^{t-1} \norm{\tbR^{t-j} \bpi_\bC} \norm{ \nabla f(\hat{x}^{(j)})}^2 .
\end{aligned}
\]
\textbf{Step 2}: For $\bH_{2,t}:= \sum_{j=0}^{t-2}\sum_{k=1}^{t-j-1}\tbR^{k}  \tbC^{t-j-k} \brk{\bG^{(j+1)} - \bG^{(j)}} + \sum_{k=1}^{t}  \tbR^{k} \tbC^{t-k}\bG^{(0)}$, we have
\[
\begin{aligned}
    & \bH_{2,t} = \sum_{j=0}^{t-2}\sum_{k=1}^{t-j-1}\tbR^{k} \tbC^{t-j-k} \brk{\btheta^{(j+1)} - \btheta^{(j)}} +\sum_{k=1}^{t}  \tbR^{k} \tbC^{t-k}\btheta^{(0)}\\
    & + \sum_{j=0}^{t-2}\sum_{k=1}^{t-j-1}\tbR^{k} \tbC^{t-j-k} \brk{\nabla \bF^{(j+1)} - \nabla \bF^{(j)}} + \sum_{k=1}^{t}  \tbR^{k} \tbC^{t-k}\nabla \bF^{(0)}  .\\
\end{aligned}
\]
Define $\bH_{2,t,1}$ and $\bH_{2,t,2}$ as the terms that involve only $\btheta$ and $\nabla \bF$ respectively. We bound $\bH_{2,t,1}$ and $\bH_{2,t,2}$ separately in step 2.1 and 2.2.\\
\textbf{Step 2.1}: It follows from the transformation in Equation \eqref{eq:pp2_PiC} that 
\[
\begin{aligned}
    \bH_{2,t,1}:=& \sum_{j=0}^{t-2} \prt{ \sum_{k=1}^{t-j}\tbR^{k}  \tbC^{t-j-k+1} - \sum_{k=1}^{t-j-1}\tbR^{k}  \tbC^{t-j-k} } \btheta^{(j)} \\
    & + \tbR\tbC\btheta^{(t-1)}.
\end{aligned}
\]
Then, by Lemma \ref{l:co-var}, we have
\[
\begin{aligned}
    & \expect\norm{\bH_{2,t,1}}_2^2 \le  n\sigma^2 \sum_{j=0}^{t-2} \norm{ \sum_{k=1}^{t-j}\tbR^{k}  \tbC^{t-j-k+1} - \sum_{k=1}^{t-j-1}\tbR^{k}\tbC^{t-j-k} }_2^2 \\
    & \qquad \qquad + n\sigma^2 \norm{\tbR\tbC}_2^2 \le  n N_7 \sigma^2.
\end{aligned}
\]
\textbf{Step 2.2}: For $\bH_{2,t,2}$, it holds that
\[
\begin{aligned}
    & \bH_{2,t,2} :=  \bH_{2,t,2,1} + \bH_{2,t,2,2} \\
    = & \sum_{j=0}^{t-2} \prt{ \sum_{k=1}^{t-j}\tbR^{k}  \tbC^{t-j-k+1} - \sum_{k=1}^{t-j-1}\tbR^{k}  \tbC^{t-j-k} } \nabla \bar{\bF}^{(j)} \\
    & + \tbR\tbC\nabla \bar{\bF}^{(t-1)}+ \sum_{k=1}^{t} \tbR^{k}  \tbC^{t-k}\nabla \bF(\bar{\bX}^{(0)})\\
    & + \sum_{j=0}^{t-2}\sum_{k=1}^{t-j-1}\tbR^{k} \tbC^{t-j-k} \brk{\nabla \bF(\bar{\bX}^{(j+1)}) - \nabla \bF(\bar{\bX}^{(j)})} ,
\end{aligned}
\]
where $\bH_{2,t,2,1}$ and $\bH_{2,t,2,2}$ are the terms that involve only $\nabla \bar{\bF}$ and $\nabla \bF(\bar{\bX})$, respectively. Then, it holds by choosing $\alpha_{t-j} = \norm{\sum_{k=1}^{t-j}\tbR^{k}  \tbC^{t-j-k+1} - \sum_{k=1}^{t-j-1}\tbR^{k}  \tbC^{t-j-k}}_2/N_8$ in Lemma \ref{lem:sum_matrix} that
\[
\begin{aligned}
    & \expect\norm{\bH_{2,t,2,1}}_2^2 =  \expect \left\| \tbR\tbC\nabla \bar{\bF}^{(t-1)} + \right.\\
    & \left.\sum_{j=0}^{t-2} \prt{ \sum_{k=1}^{t-j}\tbR^{k}  \tbC^{t-j-k+1} - \sum_{k=1}^{t-j-1}\tbR^{k}  \tbC^{t-j-k} } \nabla \bar{\bF}^{(j)}  \right\|_F^2\\
    \le & N_8 L^2\sum_{j=0}^{t-1} \expect \norm{\sum_{k=1}^{t-j}\tbR^{k}  \tbC^{t-j-k+1} - \sum_{k=1}^{t-j-1}\tbR^{k}  \tbC^{t-j-k}}_2 \norm{\hat{\bX}^{(j)}}_F^2.\\
\end{aligned}
\]
It also holds by choosing $\alpha_{t-j} = \norm{\sum_{k=1}^{t-j-1}\tbR^{k}  \tbC^{t-j-k}}_2/N_5$ in Lemma \ref{lem:sum_matrix} that 
\[
\begin{aligned}
    & \expect\norm{\bH_{2,t,2,2}}_2^2 \le 2 \norm{\sum_{k=1}^{t} \tbR^{k}  \tbC^{t-k}}_2^2 \norm{\nabla \bF^{(0)}}_F^2 + 2N_5L^2\sum_{j=0}^{t-1} \expect \norm{\sum_{k=1}^{t-j-1}\tbR^{k}  \tbC^{t-j-k}}_2 \norm{\Delta \bar{\bX}^{(j)}  }_F^2  .\\
\end{aligned}
\]

Now, taking the squared F-norm and then full expectation on both sides of \eqref{eq:PiX-0}, summing over $t$ from $0$ to $T$, and invoking Lemma \ref{lem:sum_help}, we have
\[
\begin{aligned}
    & \sum_{t=0}^{T}\expect\norm{ \hat{\bX}^{(t)} }_F^2 \le   4\gamma^2\sum_{t=0}^{T}\expect\norm{\bH_{1,t}}_F^2  + 4\gamma^2\sum_{t=0}^{T}\expect\norm{\bH_{2,t,1}}_F^2 \\
    & + 4\gamma^2\sum_{t=0}^{T}\expect\norm{\bH_{2,t,2,1}}_F^2 + 4\gamma^2\sum_{t=0}^{T}\expect\norm{\bH_{2,t,2,2}}_F^2.\\
\end{aligned}
\]
From the previous results, we have
\[
\begin{aligned}
    & \sum_{t=0}^{T}\expect\norm{ \hat{\bX}^{(t)} }_F^2 \le  8nN_4\sigma^2 \gamma^2 (T+1)  \\
    & + 4nN_7\sigma^2\gamma^2(T+1)  + 4N_8^2 L^2\gamma^2\sum_{t=0}^{T}\expect\norm{ \hat{\bX}^{(t)} }_F^2 \\
    & + 16nN_3^2L^2 \gamma^2\sum_{t=0}^{T}\expect\norm{ \hat{\bX}^{(t)} }_F^2 + 8N_5^2L^2\gamma^2 \sum_{t=0}^{T} \expect \norm{\Delta\bar{\bX}^{(t)}  }_F^2\\
    &  + 8N_6\gamma^2 \norm{\nabla \bF^{(0)}}_F^2 + 16n^2 N_3^2\gamma^2 \sum_{t=0}^{T}\expect \norm{\nabla f(\hat{x}^{(t)})}_2^2.
\end{aligned}
\]
After implementing Lemma  \ref{l:X_diff} and rearranging the terms, we get for $\gamma \le \frac{1}{6\sqrt{n}N_1 L}$ that
\[
\begin{aligned}
& \sum_{t=0}^{T}\expect\norm{ \hat{\bX}^{(t)} }_F^2 \le 4 \gamma^2 N_8^2 L^2 \sum_{t=0}^{T} \norm{\hat{\bX}^{(t)}}_F^2  + 16 \gamma^2 N_3^2 n L^2 \sum_{t=0}^{T} \expect\norm{\hat{\bX}^{(t)} }_F^2 + 4 \gamma^2 N_7 n \sigma^2\prt{T+1} \\
& + 16 \gamma^2 N_3^2 n^2 \sum_{t=0}^{T} \expect\norm{\nabla f(\hat{x}^{(t)})}_2^2 +  8 \gamma^2 N_6 \norm{\nabla \bF^{(0)}}_F^2 +  8 \gamma^2 N_4 n \sigma^2 \prt{T+1}+ 32 \gamma^4 n M_1 N_5^2 L^2 \sigma^2 (T+1) \\
&  + 800 \gamma^4 n \max\crk{ N_1^2 , n \pi^2 }N_5^2 L^4 \sum_{t=0}^{T}\expect\norm{ \hat{\bX}^{(t)} }_F^2   + 200 \gamma^4 \pi^2 n^3 N_5^2 L^2 \sum_{t=0}^{T} \expect\norm{\nabla f(\hat{x}^{(t)})}_2^2 \\
& + 100 \gamma^4 n N_2 N_5^2 L^2 \norm{\nabla \bF^{(0)}}_F^2 . \\
\end{aligned}
\]
By choosing $\gamma \le \frac{1}{20\sqrt{C_1} L}$, where $C_1 = \max\crk{nN_1^2, nN_3^2, N_8^2, \sqrt{n}N_1N_5, n\pi N_5}$, there holds
\[
\begin{aligned}
&   800 \gamma^4 n \max\crk{ N_1^2 , n \pi^2 }N_5^2 L^4  + 16 \gamma^2 N_3^2 n L^2  \\
&\  + 4 \gamma^2 N_8^2 L^2 \le \frac{5}{6} , \quad 200 \gamma^4  n^3\pi^2 N_5^2 L^2  \le 6\gamma^2n^2 \pi N_5 ,\\
& 100 \gamma^4 n N_2  N_5^2 L^2 \le 3\gamma^2 \frac{N_2 N_5}{\pi}.
\end{aligned}
\]
Hence, we have 
\[
\begin{aligned}
    & \sum_{t=0}^{T}\expect\norm{ \hat{\bX}^{(t)} }_F^2 \le 144\gamma^2 \max\crk{N_4,N_7} n \sigma^2 \prt{T+1}\\
    &\quad  + 384\gamma^4 nM_1N_5^2L^2 \sigma^2 (T+1) + \frac{5}{6}  \sum_{t=0}^{T} \norm{\hat{\bX}^{(t)}}_F^2 \\
    &\quad +  35 \gamma^2 \max\crk{ \frac{N_2N_5}{\pi}, N_6 } \norm{\nabla \bF^{(0)}}_F^2 \\
    & \quad  + 38 \gamma^2 n \max\crk{n N_3^2, n\pi N_5 } \sum_{t=0}^{T} \expect\norm{\nabla f(\hat{x}^{(t)})}_2^2.
\end{aligned}
\]
The proof is completed by rearranging the terms.

\subsection{Proof of Lemma \ref{l:descent}}
\label{pf:l:descent}
By Assumption \ref{a.smooth}, the function $f:= \frac{1}{n}\sum_{i=1}^{n} f_i$ is $L$-smooth. Then, it holds by the descent lemma that
\begin{equation}
    \label{eq:descent_lemma}
    \begin{aligned}
        &\expect f(\hat{x}^{(t+1)}) \le \expect f(\hat{x}^{(t)}) + \expect\left\langle \nabla f(\hat{x}^{(t)}), \hat{x}^{(t+1)} - \hat{x}^{(t)} \right\rangle \\
        & \qquad +\frac{L}{2} \expect \norm{\hat{x}^{(t+1)} - \hat{x}^{(t)}}^2  .
    \end{aligned}
\end{equation}
Note that, for the last term on the right hand side of Equation \eqref{eq:descent_lemma}, we have
$
\expect \norm{\hat{x}^{(t+1)} - \hat{x}^{(t)}}^2  = \frac{1}{n}\expect\norm{\Delta\bar{\bX}^{(t)}  }_F^2.
$
For the second last term, there holds
\begin{equation}
    \label{eq:second_last}
    \begin{aligned}
        & \expect\left\langle \nabla f(\hat{x}^{(t)}), \hat{x}^{(t+1)} - \hat{x}^{(t)} \right\rangle  = \expect \left\langle \nabla f(\hat{x}^{(t)}), -\gamma \bpi_{\bR}^\T\bY^{(t)} \right\rangle \\
        = & \expect \left\langle \nabla f(\hat{x}^{(t)}), -\gamma \prt{\bpi_{\bR} - \pi \mone}^\T  \bY^{(t)} \right\rangle \\
        & + \expect \left\langle \nabla f(\hat{x}^{(t)}), -\gamma \pi \mone^\T\bY^{(t)} \right\rangle.
    \end{aligned}
\end{equation}
We deal with the two terms on the right hand side of Equation \eqref{eq:second_last} separately in steps 1 and 2 below.

\textbf{Step 1}: Notice that $\prt{\bpi_{\bR} - \bpi_\bR^\T\bpi_{\bC}\mone}^\T  = \bpi_\bR^\T\bPi_{\bC}$. Then, by implementing Equation \eqref{eq:pp2_PiC}, we have
\[
\begin{aligned}
    & \expect \left\langle \nabla f(\hat{x}^{(t)}), -\gamma \prt{\bpi_{\bR} - \pi \mone}^\T \bY^{(t)} \right\rangle \\
    = &\sum_{j=1}^{t-1} \expect \left\langle \nabla f(\hat{x}^{(t)}), -\gamma \bpi_{\bR}^\T  \prt{\bC^{t-j + 1} - \bC^{t-j}}\btheta^{(j)} \right\rangle \\
    &  + \expect \left\langle \nabla f(\hat{x}^{(t)}), -\gamma \bpi_{\bR}^\T \tbC\btheta^{(t)}\right\rangle \\
    &  + \gamma \expect \left\langle \nabla f(\hat{x}^{(t)}), \sum_{j=0}^{t-1} \bpi_\bR^\T \tbC^{t-j} \brk{\nabla\bF^{(j+1)} - \nabla\bF^{(j)} } \right\rangle \\
    & + \gamma \expect \ip{\nabla f(\hat{x}^{(t)})}{ \bpi_\bR^\T \tbC^{t}\nabla \bF^{(0)}  }.
\end{aligned}
\]
We bound the above four terms in the equation step by step, as outlined in Steps 1.1 to 1.3.\\
\textbf{Step 1.1}: Under Assumption \ref{a.var},
\begin{equation}
    \begin{aligned}
        &\expect \left\langle \nabla f(\hat{x}^{(j)}), -\gamma \bpi_{\bR}^\T  \prt{\bC^{t-j+1} - \bC^{t-j}}\btheta^{(j)} \right\rangle \bigg| \cF_j \\
        =& \left\langle \nabla f(\hat{x}^{(j)}), -\expect \gamma \bpi_{\bR}^\T \prt{\bC^{t-j+1} - \bC^{t-j}}\btheta^{(j)} \bigg| \cF_j\right\rangle =0.
    \end{aligned}
    \label{eq:delay_independent}
\end{equation}
Then, we have
    \[
    \begin{aligned}
        &\expect \left\langle \nabla f(\hat{x}^{(t)}), -\gamma \bpi_{\bR}^\T \prt{\bC^{t-j+1} - \bC^{t-j}}\btheta^{(j)} \right\rangle\\
        = &  \expect \left\langle \nabla f(\hat{x}^{(t)}) - \nabla f(\hat{x}^{(j)}), -\gamma \bpi_{\bR}^\T  \prt{\bC^{t-j+1} - \bC^{t-j}}\btheta^{(j)} \right\rangle \\
        \le & \alpha \expect\norm{\nabla f(\hat{x}^{(t)}) - \nabla f(\hat{x}^{(j)}) }^2/2 \\
        & + \gamma^2 \frac{1}{2\alpha} \expect \norm{\bpi_{\bR}^\T  \prt{\bC^{t-j+1} - \bC^{t-j}}\btheta^{(j)}}^2, \\ 
    \end{aligned}
    \]
where we use the fact $\left\langle a,b \right\rangle\le \alpha \norm{a}^2/2 + \norm{b}^2/(2\alpha)$ in the first inequality and choose $\alpha = \frac{\norm{ \bpi_{\bR}^\T\prt{\bC^{t-j+1} - \bC^{t-j}}}}{(t-j)L}$ if $\norm{ \bpi_{\bR}^\T\prt{\bC^{t-j+1} - \bC^{t-j}}} \ne 0$.
Consequently, the following inequality holds no matter $\norm{ \bpi_{\bR}^\T\prt{\bC^{t-j+1} - \bC^{t-j}}} = 0$ or not:
\[
\begin{aligned}
    &\expect \left\langle \nabla f(\hat{x}^{(t)}), -\gamma \bpi_{\bR}^\T \prt{\bC^{t-j+1} - \bC^{t-j}}\btheta^{(j)}\right\rangle\\
    \le & \frac{\norm{ \bpi_{\bR}^\T\prt{\bC^{t-j+1} - \bC^{t-j}}}}{2(t-j)} L \expect \norm{\hat{x}^{(t)} - \hat{x}^{(j)}}^2 \\
    &  + \gamma^2 (t-j) \norm{ \bpi_{\bR}^\T\prt{\bC^{t-j+1} - \bC^{t-j}}}  \sigma^2 L/2.
\end{aligned}
\]
 We deal with the above two terms one by one. First, observe that
\[
\begin{aligned}
    & \norm{\hat{x}^{(t)} - \hat{x}^{(j)}}^2 = \norm{\sum_{k=j}^{t-1} \prt{\hat{x}^{(k+1)} - \hat{x}^{(k)}}}^2 \\
    \le & (t-j)\sum_{k=j}^{t-1} \norm{\hat{x}^{(k+1)} - \hat{x}^{(k)}}^2.
\end{aligned}
\]
We have
\begin{equation}
    \label{eq:second_last_2}
    \begin{aligned}
        & \sum_{j=1}^{t-1} \frac{\norm{ \bpi_{\bR}^\T\prt{\bC^{t-j+1} - \bC^{t-j}}}}{t-j} L \expect \norm{\hat{x}^{(t)} - \hat{x}^{(j)}}^2 \\
        \le & L \sum_{j=1}^{t-1} \norm{ \bpi_{\bR}^\T\prt{\bC^{t-j+1} - \bC^{t-j}}}  \sum_{k=j}^{t-1} \expect \norm{\hat{x}^{(k+1)} - \hat{x}^{(k)}}^2.
    \end{aligned}
\end{equation}
Summing over $t$ from $1$ to $T$ on both sides of Equation \eqref{eq:second_last_2} and exchanging the order of the summation, we get
\begin{equation}
    \label{eq:l_descent_tj}
    \begin{aligned}
        &\sum_{t=0}^{T}\sum_{j=1}^{t-1} \frac{\norm{ \bpi_{\bR}^\T\prt{\bC^{t-j+1} - \bC^{t-j}}}}{t-j} L \expect \norm{\hat{x}^{(t)} - \hat{x}^{(j)}}^2 \\
        \le & L \sum_{t=0}^{T} \sum_{j=1}^{t-1}\sum_{k=j}^{t-1}  \norm{ \bpi_{\bR}^\T\prt{\bC^{t-j+1} - \bC^{t-j}}}  \expect \norm{\hat{x}^{(k+1)} - \hat{x}^{(k)}}^2\\
        = & L \sum_{t=0}^{T}\sum_{k=1}^{t-1}  \sum_{j=1}^{k}  \norm{ \bpi_{\bR}^\T\prt{\bC^{t-j+1} - \bC^{t-j}}}  \expect \norm{\hat{x}^{(k+1)} - \hat{x}^{(k)}}^2\\
        = &  L \sum_{k=1}^{T-1}  \sum_{t=k+1}^{T}\sum_{j=1}^{k}  \norm{ \bpi_{\bR}^\T\prt{\bC^{t-j+1} - \bC^{t-j}}}  \expect \norm{\hat{x}^{(k+1)} - \hat{x}^{(k)}}^2.\\
    \end{aligned}
\end{equation}
By choosing $a_{t-j} = \norm{ \bpi_{\bR}^\T\prt{\bC^{t-j+1} - \bC^{t-j}}}$ in Lemma \ref{lem:series_help1}, we have
\[
\sum_{t=k+1}^{T}\sum_{j=1}^{k}  \norm{ \bpi_{\bR}^\T\prt{\bC^{t-j+1} - \bC^{t-j}}} \le \sum_{i=1}^{T-1} i a_i \le M_2.
\]
Hence, 
\[
\begin{aligned}
    & \sum_{t=0}^{T}\sum_{j=1}^{t-1} \frac{\norm{ \bpi_{\bR}^\T\prt{\bC^{t-j+1} - \bC^{t-j}}}}{t-j} L \expect \norm{\hat{x}^{(t)} - \hat{x}^{(j)}}^2 \\
    \le & \frac{M_2 L}{n} \sum_{t=0}^{T}\expect\norm{ \Delta\bar{\bX}^{(t)} }_F^2.
\end{aligned}
\]

 Second, for the term concerning $\sigma^2$, we have
\[
\sum_{j=1}^{t-1}  \gamma^2 (t-j) \norm{ \bpi_{\bR}^\T\prt{\bC^{t-j+1} - \bC^{t-j}}}  \sigma^2 L \le \gamma^2 M_2 \sigma^2 L.
\]

 Additionally, in case that $\cG_{\bR}$ and $\cG_{\bC^\T}$ are spanning trees with a common root as in \cite{you2025distributed}, by invoking Lemma 3.4 in \cite{you2025distributed}, it holds that
\[
\expect \left\langle \nabla f(\hat{x}^{(t)}), -\gamma \bpi_{\bR}^\T \prt{\bC^{t-j+1} - \bC^{t-j}}\btheta^{(j)}\right\rangle = 0.
\]
Thus, after defining $\tilde{M}_2=0$ in such a case and $\tilde{M}_2 = M_2$ otherwise, we obtain the following inequality:
\[
\begin{aligned}
    &\sum_{t=0}^{T}\sum_{j=1}^{t-1}\expect \left\langle \nabla f(\hat{x}^{(t)}), -\gamma \bpi_{\bR}^\T \prt{\bC^{t-j+1} - \bC^{t-j}}\btheta^{(j)}\right\rangle\\
    \le & \frac{\tilde{M}_2 L}{2n} \sum_{t=0}^{T}\expect\norm{ \Delta\bar{\bX}^{(t)} }_F^2 + \gamma^2 \tilde{M}_2 \sigma^2 L\prt{T+1}/2.
\end{aligned}
\]  
\textbf{Step 1.2}: It holds that $\expect \left\langle \nabla f(\hat{x}^{(t)}), -\gamma \bpi_{\bR}^\T \tbC\btheta^{(t)}\right\rangle = 0$ by invoking Equation \eqref{eq:delay_independent}.\\
\textbf{Step 1.3}: It holds by choosing $\alpha = n\pi$ in inequality $\left\langle a,b \right\rangle\le \alpha \norm{a}^2/2 + \norm{b}^2/(2\alpha)$ that
\[
\begin{aligned}
& \expect \left\langle \nabla f(\hat{x}^{(t)}), \sum_{j=0}^{t-1} \bpi_\bR^\T \tbC^{t-j} \brk{\nabla\bF^{(j+1)} - \nabla\bF^{(j)}}  \right\rangle 
\\
\le &  \frac{1}{2n\pi}\expect\norm{  \sum_{j=0}^{t-1} \bpi_\bR^\T \tbC^{t-j} \brk{\nabla\bF^{(j+1)} - \nabla\bF^{(j)} }  }_F^2\\
& + \frac{n\pi}{2}\expect \norm{\nabla f(\hat{x}^{(t)}) }^2  ,
\end{aligned}
\]
and
\[
\begin{aligned}
    & \expect \ip{\nabla f(\hat{x}^{(t)})}{ \bpi_\bR^\T \tbC^{t}\nabla \bF^{(0)} } \le \frac{n\pi }{2}\expect \norm{\nabla f(\hat{x}^{(t)}) }^2\\
    & \quad + \frac{1}{2n\pi }  \norm{\bpi_\bR^\T \tbC^{t} }^2 \norm{\nabla \bF^{(0)}}_F^2 .\\
\end{aligned}
\]
By choosing $\alpha_i = \norm{\bpi_\bR^\T \tbC^i}/ N_1$ in Lemma \ref{lem:sum_matrix}, we have
\[
\begin{aligned}
    & \expect\norm{  \sum_{j=0}^{t-1} \bpi_\bR^\T \tbC^{t-j} \brk{\nabla\bF^{(j+1)} - \nabla\bF^{(j)} }  }_F^2\\
    \le & N_1 \sum_{j=0}^{t-1}\expect \norm{\bpi_\bR^\T \tbC^{t-j}} \norm{\nabla\bF^{(j+1)} - \nabla\bF^{(j)} }_F^2 .
\end{aligned}
\]
Summing over $t$ from $1$ to $T$ on  both sides of the above inequality and invoking Assumption \ref{a.smooth} and Lemma \ref{lem:sum_help}, we have
\begin{equation}
    \label{eq:gradient_diff}
    \begin{aligned}
        & \sum_{t=0}^{T}\expect\norm{  \sum_{j=0}^{t-1} \bpi_\bR^\T \tbC^{t-j} \brk{\nabla\bF^{(j+1)} - \nabla\bF^{(j)} }  }_F^2\\
        \le & 3N_1^2 L^2 \sum_{t=0}^{T}\expect \norm{ \Delta\bar{\bX}^{(t)} }_F^2  + 6N_1^2L^2   \sum_{t=0}^{T}\expect \norm{\hat{\bX}^{(t)}}_F^2 .
    \end{aligned}
\end{equation}
\textbf{Step 2}: From $\mone\bY^{(t)} = \mone\bG^{(t)}$ and Equation \eqref{eq:delay_independent}, we have
\[
\begin{aligned}
   & \expect \left\langle \nabla f(\hat{x}^{(t)}), -\gamma \pi \mone^\T\bY^{(t)} \right\rangle  = \expect \left\langle \nabla f(\hat{x}^{(t)}), -\gamma \pi \mone^\T\btheta^{(t)} \right\rangle \\
   & +  \expect \left\langle \nabla f(\hat{x}^{(t)}), -\gamma \pi\mone^\T\nabla \bbF^{(t)} \right\rangle - \gamma n \pi \expect\norm{\nabla f(\hat{x}^{(t)})}^2\\
   & \le  -\frac{3\gamma n \pi }{4} \expect\norm{\nabla f(\hat{x}^{(t)})}_2^2  + \gamma \pi \expect \norm{ \nabla \bbF^{(t)} }_F^2.  \\
\end{aligned}
\]

 Combining the results in the steps above, it holds by summing from $t = 0$ to $T$ in Equation \eqref{eq:descent_lemma}, and grouping related quantities that
\[
\begin{aligned}
    & -\Delta_f\le   \frac{1}{2}\gamma^2 \tilde{M}_2 \sigma^2 L \prt{T+1} -\frac{\gamma n \pi }{4} \sum_{t=0}^{T} \expect\norm{\nabla f(\hat{x}^{(t)})}_2^2 \\
    & + \prt{ \frac{3\gamma N_1^2 L^2}{n\pi } +  \max\crk{\tilde{M}_2, 1}\frac{L}{n} } \sum_{t=0}^{T}\expect\norm{ \Delta\bar{\bX}^{(t)} }_F^2 \\
    & + \frac{7\gamma \max\crk{N_1^2, n \pi^2} L^2}{ n\pi }  \sum_{t=0}^{T}\expect \norm{\hat{\bX}^{(t)}}_F^2 + \frac{\gamma N_2 }{n\pi} \norm{\nabla \bF^{(0)}}_F^2 .  
\end{aligned}
\]
Denote $C_2 = 3\gamma N_1^2 L /\pi + \max\crk{\tilde{M}_2, 1}$. Then, the coefficient of the term $\sum_{t=0}^{T}\expect\norm{ \Delta\bar{\bX}^{(t)} }_F^2$ above is $C_2 L/n$.\\
Invoking Lemma \ref{l:X_diff}, we have for $\gamma \le \frac{1}{6\sqrt{n}N_1 L}$ that
\[
\begin{aligned}
        &-\Delta_f 
    \le  \prt{4C_2 M_1 + \frac{1}{2} \tilde{M}_2} \gamma^2\sigma^2 L \prt{T+1} \\
    & + \frac{C_3 L}{n} \sum_{t=0}^{T}\expect \norm{\hat{\bX}^{(t)}}_F^2  + \prt{12\gamma^2  C_2 L + \frac{\gamma  }{n\pi } }N_2 \norm{\nabla \bF^{(0)}}_F^2 \\
    &  + \prt{ 24 \gamma^2 n^2 \pi^2 C_2 L -\frac{\gamma n \pi }{4}}  \sum_{t=0}^{T} \expect\norm{\nabla f(\hat{x}^{(t)})}_2^2. 
\end{aligned}
\]
where $C_3 := 96\gamma^2 n C_2 \max\crk{N_1^2,n\pi^2} L^2 + 7\gamma \max\crk{N_1^2, n \pi^2} L/\pi $. 

Invoking Lemma \ref{l:PiX}, we have for $\gamma \le \frac{1}{20\sqrt{C_1} L} (\le \frac{1}{6\sqrt{n}N_1 L})$, where $C_1 =  \max\crk{nN_1^2, nN_3^2, N_8^2, \sqrt{n}N_1 N_5, n\pi N_5},$
that
\begin{equation}
    \label{eq:descent_result}
    \begin{aligned}
    &-\Delta_f  \le   2400C_3\gamma^4 M_1N_5^2 L^3 \sigma^2 (T+1)  \\
    & + \brk{870 C_3\max\crk{N_4, N_7} + 4C_2 M_1 + \frac{1}{2} \tilde{M}_2} \gamma^2\sigma^2 L \prt{T+1} \\
    &  + C_4 \sum_{t=0}^{T} \expect\norm{\nabla f(\hat{x}^{(t)})}^2+ C_5 \norm{\nabla \bF^{(0)}}_F^2, 
   \end{aligned}
\end{equation}
in which
\[
\begin{aligned}
    C_4 & := 228\gamma^2 C_3 \max\crk{nN_3^2, n\pi N_5}L  + 24 \gamma^2 n^2 \pi^2 C_2 L -\frac{\gamma n \pi }{4} ,\\
    C_5 & := 210\gamma^2\max\crk{\frac{N_2 N_5}{\pi}, N_6} \frac{C_3 L}{n} + \brk{12\gamma^2  C_2 L + \frac{\gamma  }{2n\pi}} N_2.
\end{aligned}
\]

We are now ready to derive the desired result.
Define the following constants:
\[
\begin{aligned}
    P_1 = & \max\crk{nN_1^2, nN_3^2,N_8^2, \sqrt{n} N_1 N_5, n\pi N_5} ,\\
    P_2 = & \max\crk{nN_1^2, n^2\pi^2, nN_1^2 \tilde{M}_2, n^2\pi^2 \tilde{M}_2^2} ,\\
    P_3 = & \max\crk{ \frac{N_1^4}{\pi^2},\frac{N_3^4}{\pi^2}, n^2 \pi^2 },  P_4 =  \sqrt{\max\crk{P_2,P_3}}N_5,\\
    Q = &\max\crk{M_1, \tilde{M}_2, M_1 \tilde{M}_2}.
\end{aligned}
\]
Then, we have for $\gamma \le \frac{1}{500\sqrt{\max\crk{P_1,P_2,P_3,P_4}}L}(\le \frac{1}{20\sqrt{C_1} L})$ that
\[
\begin{aligned}
    C_2 &\le \frac{3N_1^2}{200\pi\sqrt{P_3}}+\max\crk{\tilde{M}_2,1}  \le 2\max\crk{\tilde{M}_2,1},\\
    C_3 & \le \gamma L \prt{\sqrt{P_2} + 7\sqrt{P_3}} \le 8\gamma L\sqrt{ \max\crk{P_2,P_3}} \le \frac{1}{20}. \\
\end{aligned}
\]
Thus, it holds that
\[
\begin{aligned}
    & C_4 \le  \gamma n \pi\left( 228\gamma^2\max\crk{\sqrt{P_2 P_3},P_3,P_4}L^2 \right. \\
    &\qquad  \left. + 48\gamma  \sqrt{P_2}L -\frac{1}{4} \right) \le  \gamma n \pi\prt{\frac{1}{40} + \frac{1}{8} - \frac{1}{4}} = -\frac{\gamma n \pi}{10}.
    \end{aligned}
\]
Back to the Equation \eqref{eq:descent_result}, we have
\[
\begin{aligned}
    C_5 \le & \frac{3\gamma}{n}\max\crk{\frac{N_2}{\pi},\frac{N_2}{n\pi^2} , \frac{N_6}{N_5} } := \frac{3\gamma P_5}{n},
\end{aligned}
\]
where $P_5 := \max\crk{N_2/\pi,N_2/(n\pi^2) , N_6/N_5} $.
Then, for $\gamma \le \frac{1}{500\sqrt{\max\crk{P_1,P_2,P_3,P_4}}L}$, it holds by rearranging the terms in Equation \eqref{eq:descent_result} that
\[
\begin{aligned}
& \frac{\gamma n \pi }{10} \sum_{t=0}^{T} \expect\norm{\nabla f(\hat{x}^{(t)})}_2^2 \le \Delta_f + 10Q \gamma^2\sigma^2 L \prt{T+1} \\
&  + 8000\sqrt{\max\crk{P_2,P_3}} \max\crk{N_4, N_7} \gamma^3 \sigma^2 L^2 \prt{T+1} \\ 
& + 20000\sqrt{\max\crk{P_2,P_3}} M_1 N_5^2\gamma^5\sigma^2 L^4\prt{T+1}  + \frac{3\gamma P_5}{n} \norm{\nabla \bF^{(0)}}_F^2.\\
\end{aligned}
\]
Dividing both sides above by $\frac{\gamma n \pi }{10}$ and $T+1$, we have
\[
\begin{aligned}
& \frac{1}{T+1}\sum_{t=0}^{T} \expect\norm{\nabla f(\hat{x}^{(t)})}_2^2 \le \frac{10\Delta_f}{\gamma n \pi \prt{T+1}}  + \frac{ 80000\sqrt{\max\crk{P_2,P_3}} \max\crk{N_4, N_7} \sigma^2 L^2  }{n \pi }\gamma^2 \\
&  + \frac{200000\sqrt{\max\crk{P_2,P_3}} M_1 N_5^2\sigma^2 L^4 }{n \pi }\gamma^4 + \frac{100Q\sigma^2 L}{n \pi } \gamma  + \frac{30 P_5}{n\pi \prt{T+1}} \frac{1}{n}\norm{\nabla \bF^{(0)}}_F^2.
\end{aligned}
\]

\subsection{Proof of Theorem \ref{thm:convergence}}
\label{pf:thm:convergence}
Invoking Lemma \ref{l:descent} and referring to 
    \[
    \begin{aligned}
        A &= \frac{\Delta_f}{n\pi} , B  = \frac{Q\sigma^2 L}{n\pi}, \alpha  = 500\sqrt{\max\crk{P_1,P_2,P_3,P_4}}L, C = \frac{\sqrt{\max\crk{P_2,P_3}} \max\crk{N_4, N_7} \sigma^2 L^2}{n\pi},\\
        \gamma & = \min\crk{\prt{\frac{A}{B(T+1)}}^{\frac{1}{2}} , \prt{\frac{A}{C(T+1)}}^{\frac{1}{3}},\frac{1}{\alpha}},\\
    \end{aligned}
    \]
    it holds that
    \[
    \begin{aligned}
        &\frac{1}{T+1}\sum_{t=0}^{T} \expect\norm{\nabla f(\hat{x}^{(t)})}_2^2 \le  C_0 \prt{\prt{\frac{AB}{T+1}}^{\frac{1}{2}} + C^{\frac{1}{3}}\prt{\frac{A}{T+1}}^{\frac{2}{3}} + \frac{\alpha A}{T+1}} \\
       &  + \frac{30 P_5}{n\pi\prt{T+1}}\frac{1}{n}\norm{\nabla \bF^{(0)}}_F^2  + C_0\frac{\sqrt{\max\crk{P_2,P_3}} M_1 N_5^2\sigma^2 L^4 }{n \pi }\prt{\frac{A}{B(T+1)}}^2,    
    \end{aligned}
    \]
    where $C_0 = 2\times 10^6$. This completes the proof.

\bibliographystyle{siamplain}
\bibliography{references_all}

\begin{thebibliography}{10}

\bibitem{alghunaim2022unified}
{\sc S.~A. Alghunaim and K.~Yuan}, {\em A unified and refined convergence analysis for non-convex decentralized learning}, IEEE Transactions on Signal Processing, 70 (2022), pp.~3264--3279.

\bibitem{assran2019stochastic}
{\sc M.~Assran, N.~Loizou, N.~Ballas, and M.~Rabbat}, {\em Stochastic gradient push for distributed deep learning}, in International Conference on Machine Learning, PMLR, 2019, pp.~344--353.

\bibitem{chen2012diffusion}
{\sc J.~Chen and A.~H. Sayed}, {\em Diffusion adaptation strategies for distributed optimization and learning over networks}, IEEE Transactions on Signal Processing, 60 (2012), pp.~4289--4305.

\bibitem{gharesifard2010does}
{\sc B.~Gharesifard and J.~Cort{\'e}s}, {\em When does a digraph admit a doubly stochastic adjacency matrix?}, in Proceedings of the 2010 American Control Conference, IEEE, 2010, pp.~2440--2445.

\bibitem{horn2012matrix}
{\sc R.~A. Horn and C.~R. Johnson}, {\em Matrix analysis}, Cambridge university press, 2012.

\bibitem{huang2024distributed}
{\sc K.~Huang, X.~Li, and S.~Pu}, {\em Distributed stochastic optimization under a general variance condition}, IEEE Transactions on Automatic Control,  (2024).

\bibitem{huang2022improving}
{\sc K.~Huang and S.~Pu}, {\em Improving the transient times for distributed stochastic gradient methods}, IEEE Transactions on Automatic Control, 68 (2022), pp.~4127--4142.

\bibitem{koloskova2021improved}
{\sc A.~Koloskova, T.~Lin, and S.~U. Stich}, {\em An improved analysis of gradient tracking for decentralized machine learning}, Advances in Neural Information Processing Systems, 34 (2021), pp.~11422--11435.

\bibitem{koloskova2020unified}
{\sc A.~Koloskova, N.~Loizou, S.~Boreiri, M.~Jaggi, and S.~Stich}, {\em A unified theory of decentralized sgd with changing topology and local updates}, in International Conference on Machine Learning, PMLR, 2020, pp.~5381--5393.

\bibitem{koloskova2019decentralized}
{\sc A.~Koloskova, S.~Stich, and M.~Jaggi}, {\em Decentralized stochastic optimization and gossip algorithms with compressed communication}, in International Conference on Machine Learning, PMLR, 2019, pp.~3478--3487.

\bibitem{kungurtsev2023decentralized}
{\sc V.~Kungurtsev, M.~Morafah, T.~Javidi, and G.~Scutari}, {\em Decentralized asynchronous nonconvex stochastic optimization on directed graphs}, IEEE Transactions on Control of Network Systems, 10 (2023), pp.~1796--1804.

\bibitem{lecun2010mnist}
{\sc Y.~LeCun, C.~Cortes, C.~Burges, et~al.}, {\em Mnist handwritten digit database}, 2010.

\bibitem{li2019decentralized}
{\sc Z.~Li, W.~Shi, and M.~Yan}, {\em A decentralized proximal-gradient method with network independent step-sizes and separated convergence rates}, IEEE Transactions on Signal Processing, 67 (2019), pp.~4494--4506.

\bibitem{lian2017can}
{\sc X.~Lian, C.~Zhang, H.~Zhang, C.-J. Hsieh, W.~Zhang, and J.~Liu}, {\em Can decentralized algorithms outperform centralized algorithms? a case study for decentralized parallel stochastic gradient descent}, Advances in neural information processing systems, 30 (2017).

\bibitem{liang2023understanding}
{\sc L.~Liang, X.~Huang, R.~Xin, and K.~Yuan}, {\em Understanding the influence of digraphs on decentralized optimization: Effective metrics, lower bound, and optimal algorithm}, arXiv preprint arXiv:2312.04928,  (2023).

\bibitem{liao2024robust}
{\sc Y.~Liao, Z.~Li, S.~Pu, and T.-H. Chang}, {\em A robust compressed push-pull method for decentralized nonconvex optimization}, arXiv preprint arXiv:2408.01727,  (2024).

\bibitem{nadiradze2021asynchronous}
{\sc G.~Nadiradze, A.~Sabour, P.~Davies, S.~Li, and D.~Alistarh}, {\em Asynchronous decentralized sgd with quantized and local updates}, Advances in Neural Information Processing Systems, 34 (2021), pp.~6829--6842.

\bibitem{nedic2014distributed}
{\sc A.~Nedi{\'c} and A.~Olshevsky}, {\em Distributed optimization over time-varying directed graphs}, IEEE Transactions on Automatic Control, 60 (2014), pp.~601--615.

\bibitem{nedic2016stochastic}
{\sc A.~Nedi{\'c} and A.~Olshevsky}, {\em Stochastic gradient-push for strongly convex functions on time-varying directed graphs}, IEEE Transactions on Automatic Control, 61 (2016), pp.~3936--3947.

\bibitem{nedic2018network}
{\sc A.~Nedi{\'c}, A.~Olshevsky, and M.~G. Rabbat}, {\em Network topology and communication-computation tradeoffs in decentralized optimization}, Proceedings of the IEEE, 106 (2018), pp.~953--976.

\bibitem{nedic2017achieving}
{\sc A.~Nedic, A.~Olshevsky, and W.~Shi}, {\em Achieving geometric convergence for distributed optimization over time-varying graphs}, SIAM Journal on Optimization, 27 (2017), pp.~2597--2633.

\bibitem{nedic2009distributed}
{\sc A.~Nedic and A.~Ozdaglar}, {\em Distributed subgradient methods for multi-agent optimization}, IEEE Transactions on Automatic Control, 54 (2009), pp.~48--61.

\bibitem{nguyen2023accelerated}
{\sc D.~T.~A. Nguyen, D.~T. Nguyen, and A.~Nedi{\'c}}, {\em Accelerated $ ab $/push--pull methods for distributed optimization over time-varying directed networks}, IEEE Transactions on Control of Network Systems, 11 (2023), pp.~1395--1407.

\bibitem{pu2021distributed}
{\sc S.~Pu and A.~Nedi{\'c}}, {\em Distributed stochastic gradient tracking methods}, Mathematical Programming, 187 (2021), pp.~409--457.

\bibitem{pu2020asymptotic}
{\sc S.~Pu, A.~Olshevsky, and I.~C. Paschalidis}, {\em Asymptotic network independence in distributed stochastic optimization for machine learning: Examining distributed and centralized stochastic gradient descent}, IEEE signal processing magazine, 37 (2020), pp.~114--122.

\bibitem{pu2021sharp}
{\sc S.~Pu, A.~Olshevsky, and I.~C. Paschalidis}, {\em A sharp estimate on the transient time of distributed stochastic gradient descent}, IEEE Transactions on Automatic Control, 67 (2021), pp.~5900--5915.

\bibitem{pu2020push}
{\sc S.~Pu, W.~Shi, J.~Xu, and A.~Nedi{\'c}}, {\em Push--pull gradient methods for distributed optimization in networks}, IEEE Transactions on Automatic Control, 66 (2020), pp.~1--16.

\bibitem{ren2005consensus}
{\sc W.~Ren and R.~W. Beard}, {\em Consensus seeking in multiagent systems under dynamically changing interaction topologies}, IEEE Transactions on automatic control, 50 (2005), pp.~655--661.

\bibitem{saadatniaki2020decentralized}
{\sc F.~Saadatniaki, R.~Xin, and U.~A. Khan}, {\em Decentralized optimization over time-varying directed graphs with row and column-stochastic matrices}, IEEE Transactions on Automatic Control, 65 (2020), pp.~4769--4780.

\bibitem{scutari2019distributed}
{\sc G.~Scutari and Y.~Sun}, {\em Distributed nonconvex constrained optimization over time-varying digraphs}, Mathematical Programming, 176 (2019), pp.~497--544.

\bibitem{shi2015extra}
{\sc W.~Shi, Q.~Ling, G.~Wu, and W.~Yin}, {\em Extra: An exact first-order algorithm for decentralized consensus optimization}, SIAM Journal on Optimization, 25 (2015), pp.~944--966.

\bibitem{song2022communication}
{\sc Z.~Song, W.~Li, K.~Jin, L.~Shi, M.~Yan, W.~Yin, and K.~Yuan}, {\em Communication-efficient topologies for decentralized learning with $ o (1) $ consensus rate}, Advances in Neural Information Processing Systems, 35 (2022), pp.~1073--1085.

\bibitem{spiridonoff2020robust}
{\sc A.~Spiridonoff, A.~Olshevsky, and I.~C. Paschalidis}, {\em Robust asynchronous stochastic gradient-push: Asymptotically optimal and network-independent performance for strongly convex functions}, Journal of machine learning research, 21 (2020), pp.~1--47.

\bibitem{tang2018d}
{\sc H.~Tang, X.~Lian, M.~Yan, C.~Zhang, and J.~Liu}, {\em $ d2$: Decentralized training over decentralized data}, in International Conference on Machine Learning, PMLR, 2018, pp.~4848--4856.

\bibitem{tian2020achieving}
{\sc Y.~Tian, Y.~Sun, and G.~Scutari}, {\em Achieving linear convergence in distributed asynchronous multiagent optimization}, IEEE Transactions on Automatic Control, 65 (2020), pp.~5264--5279.

\bibitem{xi2017add}
{\sc C.~Xi, R.~Xin, and U.~A. Khan}, {\em Add-opt: Accelerated distributed directed optimization}, IEEE Transactions on Automatic Control, 63 (2017), pp.~1329--1339.

\bibitem{xin2018linear}
{\sc R.~Xin and U.~A. Khan}, {\em A linear algorithm for optimization over directed graphs with geometric convergence}, IEEE Control Systems Letters, 2 (2018), pp.~315--320.

\bibitem{xin2020general}
{\sc R.~Xin, S.~Pu, A.~Nedi{\'c}, and U.~A. Khan}, {\em A general framework for decentralized optimization with first-order methods}, Proceedings of the IEEE, 108 (2020), pp.~1869--1889.

\bibitem{xin2019distributed}
{\sc R.~Xin, A.~K. Sahu, U.~A. Khan, and S.~Kar}, {\em Distributed stochastic optimization with gradient tracking over strongly-connected networks}, in 2019 IEEE 58th Conference on Decision and Control (CDC), IEEE, 2019, pp.~8353--8358.

\bibitem{yang2019survey}
{\sc T.~Yang, X.~Yi, J.~Wu, Y.~Yuan, D.~Wu, Z.~Meng, Y.~Hong, H.~Wang, Z.~Lin, and K.~H. Johansson}, {\em A survey of distributed optimization}, Annual Reviews in Control, 47 (2019), pp.~278--305.

\bibitem{ying2021exponential}
{\sc B.~Ying, K.~Yuan, Y.~Chen, H.~Hu, P.~Pan, and W.~Yin}, {\em Exponential graph is provably efficient for decentralized deep training}, Advances in Neural Information Processing Systems, 34 (2021), pp.~13975--13987.

\bibitem{you2024b}
{\sc R.~You and S.~Pu}, {\em B-ary tree push-pull method is provably efficient for decentralized learning on heterogeneous data}, arXiv preprint arXiv:2404.05454,  (2024).

\bibitem{you2025distributed}
{\sc R.~You and S.~Pu}, {\em Distributed learning over arbitrary topology: Linear speed-up with polynomial transient time}, arXiv preprint arXiv:2503.16123,  (2025).

\bibitem{yuan2023removing}
{\sc K.~Yuan, S.~A. Alghunaim, and X.~Huang}, {\em Removing data heterogeneity influence enhances network topology dependence of decentralized sgd}, Journal of Machine Learning Research, 24 (2023), pp.~1--53.

\bibitem{zhao2023asymptotic}
{\sc S.~Zhao and Y.~Liu}, {\em Asymptotic properties of s-ab method with diminishing step-size}, IEEE Transactions on Automatic Control,  (2023).

\bibitem{zhu2024r}
{\sc Z.~Zhu, Y.~Tian, Y.~Huang, J.~Xu, and S.~He}, {\em R-fast: Robust fully-asynchronous stochastic gradient tracking over general topology}, IEEE Transactions on Signal and Information Processing over Networks,  (2024).

\end{thebibliography}

\end{document}